\input ./style/arxiv-general.cfg
\documentclass[aap,MSNbibl,nameyear,seceqn,dvips]{arximspdf}
\makeatletter
   \@ifpackageloaded{graphicx}{}{\usepackage{graphicx}}
\makeatother
\usepackage{algorithm}
\usepackage{algorithmic}

\doi{10.1214/15-AAP1105}
\volume{26}
\issue{2}
\pubyear{2016}
\firstpage{818}
\lastpage{856}
\docsubty{FLA}

\makeatletter

\def\afrac#1#2{#1/(#2)}

\newcommand{\rrvert}{\vert}
\newcommand{\rrVert}{\Vert}
\newcommand{\llvert}{\vert}
\newcommand{\llVert}{\Vert}
\newtheorem{teo}{Theorem}
\newtheorem{lemma}{Lemma}
\newtheorem{prop}{Proposition}
\newproclaim{defn}{Definition}
\newproclaim{rem}{Remark}
\newtheorem{cor}{Corollary}
\renewcommand{\mid}{|}
\def\E{\mathbf{E}}
\def\N{\mathbf{N}}
\def\R{\mathbf{R}}
\renewcommand{\E}{E}
\makeatother

\begin{document}
\begin{frontmatter}

\title{A weak approximation with asymptotic expansion and multidimensional Malliavin weights\thanksref{T1}}
\runtitle{Weak approximation with asymptotic expansion}

\begin{aug}
\author[A]{\fnms{Akihiko}~\snm{Takahashi}\ead[label=e1]{akihikot@e.u-tokyo.ac.jp}}
\and
\author[B]{\fnms{Toshihiro}~\snm{Yamada}\corref{}\ead[label=e2]{toshihiro.ymd@gmail.com}}
\runauthor{A. Takahashi and T. Yamada}
\affiliation{University of Tokyo, University of Tokyo and MTEC}
\address[A]{Graduate School of Economics\\
University of Tokyo\\
7-3-1 Hongo, Bunkyo, Tokyo 113-0033\\
Japan\\
\printead{e1}}
\address[B]{Graduate School of Economics\\
University of Tokyo\\
7-3-1 Hongo, Bunkyo, Tokyo 113-0033\\
Japan\\
and\\
Mitsubishi UFJ Trust Investment\\
\quad Technology Institute Co., Ltd. (MTEC)\\
2-6-4, Akasaka, Minato, Tokyo 107-0052\\
Japan\\
\printead{e2}}
\end{aug}
\thankstext{T1}{Supported by JSPS KAKENHI Grant Number 25380389, CARF
(Center for Advanced Research in Finance), Graduate School of
Economics, University of Tokyo and grants for Excellent Graduate
Schools, MEXT, Japan.}

%
\received{\smonth{12} \syear{2013}}
%
\revised{\smonth{9} \syear{2014}}

%
\begin{abstract}
This paper develops a new efficient scheme for approximations
of expectations of the solutions to stochastic differential equations (SDEs).
In particular, we present a method for connecting approximate operators
based on an asymptotic expansion with multidimensional Malliavin weights
to compute a target expectation value precisely.
The mathematical validity is given based on Watanabe and Kusuoka
theories in Malliavin calculus.
Moreover, numerical experiments for option pricing under local and
stochastic volatility
models confirm the effectiveness of our scheme. Especially, our weak
approximation substantially improves the accuracy at deep
Out-of-The-Moneys (OTMs).
\end{abstract}

%
\end{frontmatter}

\section{Introduction}\label{sec1}
Developing an approximation method for expectations of diffusion processes
is an interesting topic in various research fields.
In fact, it seems so useful that a precise
approximation for the expectation
would lead to substantial reduction of
computational burden so that the subsequent analyses could be very
easily implemented.
Particularly, in finance
it has drawn much attention
for
more than the past two decades
since fast and precise computation
is so important in terms of competition and risk management
in practice such as in trading and investment.

An example among a large number of the related researches
is an asymptotic expansion approach,
which is mathematically justified by
Watanabe theory [\citet{Watanabe}]
in Malliavin calculus [e.g., \citet{Malliavin}].
Especially,
the asymptotic expansion have been applied to a broad class of problems
in finance;
for instance, see
\citeauthor{Takahashi-Yamadaa} (\citeyear{Takahashi-Yamadaa,Takahashi-Yamadab,Takahashi-Yamadac,Takahashi-Yamadad})
and references therein.

Although the asymptotic expansion up to the fifth order
is known to be sufficiently accurate for option pricing
[e.g., \citet{Takahashi-Takehara-Todab}],
the main criticism against the method
would be that the approximate density function
deviates from the true density
at its tails that is, some region of the very deep Out-of-The-Money (OTM).
However, there exist similar problems, at least implicitly
in other well-known approximation methods such as \citet{Haganetal}.

On the other hand, the Monte Carlo simulation method is quite popular
mainly due to the ease of its implementation.
Nevertheless, in order to achieve accuracy sufficient enough in practice,
there exists an unavoidable drawback in computational cost
under the standard weak approximation schemes of SDEs such as the
Euler--Maruyama scheme.

To overcome this problem, \citeauthor{Kusuokaa} (\citeyear{Kusuokaa,Kusuokac,Kusuokad}) developed a high
order weak approximation scheme for SDEs
based on Malliavin calculus and
Lie algebra, which
opened the door for the possibility that the computational speed and
the accuracy
in the Monte Carlo simulation satisfies stringent requirements in
financial business.
Independently,
\citet{Lyons-Victoir} developed a cubature method on the Wiener space.
Since then, there have been a large number of researches
for weak approximations and its applications to the computational
finance inspired by those pioneering works.
For instance, see \citet{Crisan-Manolarakis-Nee} for the Kusuoka's method and its
related works [e.g., \citet{Bayer-Friz-Loeffen}].

This paper develops a new weak approximation scheme for
expectations of functions of the solutions to SDEs.
In particular,
the scheme connects approximate operators constructed based on the
asymptotic expansion.
More concretely,
a diffusion semigroup is defined as the expectation of an appropriate
function of
the solution to a certain SDE, for example, $P^{\varepsilon}_tf(x)=E[
f(X_t^{x,\varepsilon})]$
with the solution $X^{x,\varepsilon}_t$ of a SDE with perturbation
parameter $\varepsilon$ and a function $f$.
Then we approximate $P^{\varepsilon}_{t}$ by an operator
$Q^{\varepsilon,m}_{t}$
which is constructed based on the asymptotic expansion up to a certain
order $m$.
Thus, given a partition of $[0,T]$, $\pi=\{ (t_0,t_1,\ldots,t_n)\dvtx
0=t_0<t_1<\cdots<t_n=T\}$,
we are able to approximate $P^{\varepsilon}_Tf(x)$
by connecting the expansion-based approximations sequentially, that is,
with $s_k = t_k-t_{k-1}$, $k=1,\ldots,n$,
\begin{eqnarray*}
P^{\varepsilon}_Tf(x) &\simeq& Q^{\varepsilon,m}_{s_n}
Q^{\varepsilon,m}_{s_{n-1}}\cdots Q^{\varepsilon,m}_{s_1} f(x).
\end{eqnarray*}
This paper justifies this idea by applying Malliavin calculus,
particularly, theories developed by \citet{Watanabe}
and \citeauthor{Kusuokaa} (\citeyear{Kusuokaa,Kusuokab,Kusuokad}).

Moreover, we show through numerical examples for option pricing that
very few partitions such as $n=2$ is mostly enough
to substantially improve the errors at deep OTMs of expansions with
order $m=1, 2$.
For a related but different approach with similar motivation, see
Section~5 in \citet{Fujii}.

The organization of the paper is as follows. The next section
introduces the setup and
the basic results necessary for the subsequent analysis.
Section~\ref{sec3} shows our main result for a new weak approximation of the
expectation of
diffusion processes. Section~\ref{sec4} briefly describes an example for the
implementation method of our scheme,
Section~\ref{sec5} provides numerical experiments for option pricing under local
and stochastic volatility models.
Section~\ref{sec6} makes concluding remarks. The \hyperref[appe]{Appendix} gives the proofs of
Theorems \ref{11111111111111}, \ref{222222222222222222} and \ref{33333333333333} as well as Lemma~\ref{le2} and its proof.
\section{Preparation}\label{sec2}
Let $({\mathcal W},H,\Bbb{P})$ be the $d$-dimensional Wiener space,
{that is}, ${\mathcal W}=\{w\in C([0,T]\rightarrow{\mathbf R}^d);
w(0)=0 \}$ which is a real Banach space under the supremum norm, $H=\{h
\in{\mathcal W}; t\mapsto h(t)$ is absolutely  continuous and  $\llVert   h\rrVert  _H^2=\int_0^T \llvert  \frac{d}{dt}h(t)\rrvert  ^2\,dt < \infty\}$ is
a real Hilbert space under $\llVert   \cdot\rrVert  _H$ called the Cameron--Martin
subspace and $\Bbb{P}$ is the $d$-dimensional Wiener measure. Let
$B_t=(B_t^1,\ldots,B_t^d)^{\top}$ be a $d$-dimensional Brownian motion.
In this paper, we consider the following general perturbed
$N$-dimensional stochastic differential equation with $\varepsilon\in(0,1]$:
%
\begin{eqnarray}
X_{t}^{x,\varepsilon}&=&x + \int_{0}^{t}V_{0}
\bigl(\varepsilon, X_{s}^{x,\varepsilon}\bigr) \,ds +\varepsilon\sum
_{j=1}^{d} \int_{0}^{t}V_{j}
\bigl(X_{s}^{x,\varepsilon}\bigr) \,dB_{s}^{j},
\end{eqnarray}
where $V_0 \in C_b^{\infty}((0,1]\times\R^N; \R^N)$ and $V_j \in
C_b^{\infty}(\R^N; \R^N)$, $j=1,\ldots,d$ are bounded. Hereafter,
we will use the notation $Vf(x)=\sum_{i=1}^N V^i(x)(\partial
f/\partial x_i)(x)$ for $V \in C_b^{\infty}(\R^N; \R^N)$ and $f$ a
differentiable function $\R^N$ into $\R$.
$X_t^{x,\varepsilon}$ can be written in the Stratonovich form:
%
\begin{eqnarray}
X_{t}^{x,\varepsilon}&=&x+ \int_{0}^{t}
\tilde{V}_{0}\bigl(\varepsilon, X_{s}^{x,\varepsilon}\bigr)
\,ds + \varepsilon\sum_{j=1}^{d}\int
_{0}^{t}\tilde{V}_{j}
\bigl(X_{s}^{x,\varepsilon}\bigr) \circ dB_{s}^{j},
\end{eqnarray}
where
%
\begin{eqnarray}
\tilde{V}_0^i(\varepsilon,x)&=& V_0^i(
\varepsilon,x)-\frac
{\varepsilon^2}{2}\sum_{j=1}^d
V_j V_j^i(x),
\\
\varepsilon\tilde{V}_j^i(x)&=&\varepsilon
V_j^i(x),\qquad j=1,\ldots, d.
\end{eqnarray}
Here, we consider the case $V_0^i(\varepsilon,x)=\varepsilon^k \hat
{V}_0^i(x)$, $\hat{V}_0 \in C_b^{\infty}(\R^N; \R^N)$,
$k=0,1,2$, for $i=1,\ldots,N$, which is useful in applications [see
\citet{Takahashi-Toda} for the details].
Moreover, we assume the following condition~\textup{[\textbf{H}]} on the vector
fields, which ensures both the integration by parts on the Wiener space
and the asymptotic expansion in the next section.

\begin{longlist}
\item[{\textup{[\textbf{H}]}}] The matrix $A(x)=(A^{i,i'}(x))_{i,i'}$ defined by
%
\begin{eqnarray}
A^{i,i'}(x)&=&\sum_{j=1}^d
V_j^i(x)V_j^{i'}(x)\qquad\mbox{for all } x
\in{\mathbf R}^N, 1\leq i,i' \leq N
\end{eqnarray}
is nondegenerate, {that is}, $\det(A(x))>0$.
\end{longlist}

\subsection{The space \texorpdfstring{$\mathcal{K}_r$}{$\mathcal{K}_r$}}\label{sec2.1}
Let ${\mathbf D}^{k,p}(E)$, $k\geq1$, $p \in[1,\infty)$ be the space of
$k$-times Malliavin differentiable Wiener functionals $F \in
L^p({\mathcal W},E)$, where $E$ is a separable Hilbert space. See
\citet{Watanabe}, \citet{Ikeda-Watanabe}, \citet{Malliavin}, \citet{Malliavin-Thalmaier} and \citet{Nualart} for more details of the notation.
This subsection introduces the space of Wiener functionals ${\mathcal
K}_r$ developed by \citet{Kusuokab} and its properties.
The element of ${\mathcal K}_r$ is called the \textit{Kusuoka--Stroock function}.
See \citeauthor{Neea} (\citeyear{Neea,Neeb})
and \citet{Crisan-Manolarakis-Nee} for more details of the notation
and the proofs.

\begin{defn}
Given $r\in{\mathbf R}$ and $n \in{\mathbf N}$, we denote by ${\mathcal
K}_r(E,n)$ the set of functions $G\dvtx  (0,1] \times{\mathbf R}^N \rightarrow
{\mathbf D}^{n,\infty}(E)$ satisfying the following:
\begin{longlist}
\item[1.] $G(t,\cdot)$ is $n$-times continuously differentiable and
$[\partial^{\alpha}G /\partial x^{\alpha} ]$ is continuous in $(t,x)
\in(0,1] \times{\mathbf R}^N$ a.s. for any multiindex $\alpha=\alpha
^{(l)}\in\{1,\ldots, d \}^l$
with length $\llvert  \alpha\rrvert  =l \leq n$. Here, $[\partial^{\alpha}G
/\partial x^{\alpha} ]$ is the partial derivative of $G(t,x)$ given by
$\frac{\partial^l }{ \partial x_{\alpha_1} \cdots\partial x_{\alpha
_l} }G(t,x)$.\vspace*{1pt}

\item[2.] For all $k \leq n-\llvert  \alpha\rrvert  $, $p \in[1,\infty)$,
%
\begin{equation}
\sup_{t \in(0,1], x \in{\mathbf R}^{N} } t^{-r/2} \biggl\llVert \frac
{\partial^{\alpha}G }{\partial x^{\alpha}}
(t,x) \biggr\rrVert _{{\mathbf
D}^{k,p}} <\infty.
\end{equation}
\end{longlist}
We write ${\mathcal K}_r$ for ${\mathcal K}_r({\mathbf R},\infty)$.
\end{defn}

Next, we show the basic properties of the Kusuoka--Stroock functions.

\begin{lemma}[(Properties of Kusuoka--Stroock functions)]\label{Properties}\label{le1}
\hspace*{-2pt}1. The function $(t,x)\in(0,1] \times{\mathbf R}^N \mapsto
X_t^{x,\varepsilon}$ belongs to ${\mathcal K}_0$.

2. Suppose $G \in{\mathcal K}_r(n)$ where $r \geq0$. Then, for
$i=1,\ldots, d$,
%
\begin{eqnarray}\label{property1}
&& \mbox{\textup{(a)}}\quad \int_0^{\cdot}G(s,x)\,dB_s^{i}
\in{\mathcal K}_{r+1}(n)\quad\mbox{and}
\nonumber\\[-8pt]\\[-8pt]\nonumber
&&\mbox{\textup{(b)}}\quad \int_0^{\cdot}G(s,x)\,ds
\in{\mathcal K}_{r+2}(n).
\end{eqnarray}

3. If $G_i \in{\mathcal K}_{r_i}(n_i)$, $i=1,\ldots, l$, then
\begin{eqnarray}\label{property2}
&&\mbox{\textup{(a)}}\quad \prod_i^{l} G_i \in{
\mathcal K}_{r_1+\cdots+r_l}\Bigl(\min_i n_i
\Bigr)\quad\mbox{and}
\nonumber\\[-8pt]\\[-8pt]\nonumber
&&\mbox{\textup{(b)}}\quad  \sum_{i=1}^{l}G_i
\in{\mathcal K}_{\min_i r_i} \Bigl(\min_i n_i
\Bigr).
\end{eqnarray}
\end{lemma}

Then we summarize the Malliavin's integration by parts formula using
Ku\-suoka--Stroock functions. Hereafter,
for any multiindex $\alpha=\alpha^{(k)}:=(\alpha_1,\ldots,\alpha
_k) \in\{1,\ldots,N \}^k$, $k\geq1$ with the length $\llvert  \alpha^{(k)}\rrvert  =k$,
we denote by $\partial_{\alpha^{(k)}}$ the partial derivative $\frac
{\partial^{\alpha}}{\partial x^{\alpha}}=\frac{\partial^k}{
\partial x_{\alpha_1} \cdots\partial x_{\alpha_k} }$.

\begin{prop}\label{IBPKusuoka}\label{pr1}
Suppose that condition \textup{[\textbf{H}]} holds.
Let $G\dvtx  (0,1] \times{\mathbf R}^N \rightarrow{\mathbf D}^{\infty}={\mathbf
D}^{\infty,\infty}({\mathbf R})$ be an\vspace*{1pt} element of ${\mathcal K}_r$
and let $f$ be a function that belongs to the space $C_{b}^{\infty
}({\mathbf R}^{N};\R)$.
Then\vspace*{1pt} for any multiindex $\alpha^{(k)} \in\{1,\ldots,N \}^k$, $k\geq
1$, there exists $H_{\alpha^{(k)}} (X_t^{x,\varepsilon},G(t,x)) \in
{\mathcal K}_{r-k}$ such that
%
\begin{eqnarray}\label{IBP}
&&E \bigl[ {\partial_{\alpha^{(k)}}} f\bigl(X_t^{x,\varepsilon}
\bigr)G(t,x) \bigr]=E \bigl[f\bigl(X_t^{x,\varepsilon}\bigr)
H_{\alpha^{(k)}} \bigl(X_t^{x,\varepsilon},G(t,x)\bigr) \bigr],
\nonumber\\[-8pt]\\[-8pt]
\eqntext{t \in(0,1],}
\end{eqnarray}
with
%
\begin{equation}
\sup_{x \in{\mathbf R}^N } \bigl\llVert H_{\alpha^{(k)}}
\bigl(X_t^{x,\varepsilon},G(t,x)\bigr) \bigr\rrVert _{L^p}
\leq t^{(r-k)/2}{C},
\end{equation}
where $H_{\alpha^{(k)}} (X_t^{x,\varepsilon},G(t,x))$ is recursively
given by
%
\begin{eqnarray}
H_{(i)} \bigl(X_t^{x,\varepsilon},G(t,x)\bigr)&=&\delta
\Biggl(\sum_{j=1}^{N}G(t,x)
\gamma_{ij}^{X_t^{x,\varepsilon}}DX_t^{x,\varepsilon,j} \Biggr),
\\
H_{\alpha^{(k)}} \bigl(X_t^{x,\varepsilon},G(t,x)
\bigr)&=&H_{(\alpha_k)} \bigl(X_t^{x,\varepsilon},
H_{\alpha^{(k-1)}}\bigl(X_t^{x,\varepsilon},G(t,x)\bigr)\bigr),
\end{eqnarray}
and a positive constant $C$.
Here, $\delta$ is the Skorohod integral and $DX_t^{x,\varepsilon}$ is
the Malliavin derivative of $X_t^{x,\varepsilon}$,
\begin{eqnarray}
\bigl\langle DX_t^{x,\varepsilon},h \bigr\rangle_{H}&=&
\sum_{k=1}^d \int_0^t
D_{s,k}X_t^{x,\varepsilon} \frac{d}{ds}h_k(s)\,ds
\nonumber\\[-8pt]\\[-8pt]\nonumber
&=& \lim_{\lambda
\rightarrow0}\frac{X_t^{x,\varepsilon}(w+\lambda
h)-X_t^{x,\varepsilon}(w)}{\lambda},\qquad h \in H,
\nonumber
\end{eqnarray}
and $\gamma^{X_t^{x,\varepsilon}}=(\gamma_{ij}^{X_t^{x,\varepsilon
}})_{1\leq i,j \leq N}$ is the inverse matrix of the Malliavin
covariance of~$X_t^{x,\varepsilon}$.
\end{prop}

\begin{pf}
By 1, 2, 3 of Lemma \ref{Properties}, we can see that the Malliavin
covariance of $X_t^{x,\varepsilon}$ is given by
%
\begin{equation}
\sigma^{X_t^{x,\varepsilon}}_{i,j}=\sum_{k=1}^d
\int_0^t D_{s,k}X_t^{x,\varepsilon,i}D_{s,k}X_t^{x,\varepsilon,j}\,ds
\in {\mathcal K}_2,
\end{equation}
since $D_{s,k}X_t^{x,\varepsilon,i} \in{\mathcal K}_0$, $s \leq t$,
$k=1,\ldots,d$, $i=1,\ldots, N$.
Under \textup{[\textbf{H}]}, it can be shown that the nondegenerate condition of
the Malliavin covariance matrix is satisfied when $\varepsilon>0$ (but
not satisfied when $\varepsilon=0$, that is, the Malliavin covariance
matrix $\sigma^{X_t^{x,\varepsilon}}$ is not uniformly nondegenerate
in $\varepsilon$) and then (\ref{IBP}) holds [see the proofs of
Proposition 5.8, Theorems 5.9~and~6.7 of \citet{Sgigekawa}].
Also, we have $\gamma^{X_t^{x,\varepsilon}} \in{\mathcal K}_{-2}$
since $\gamma^{X_t^{x,\varepsilon}}=(\sigma^{X_t^{x,\varepsilon
}})^{-1}=\frac{\operatorname{adj} \sigma^{X_t^{x,\varepsilon}} }{\det\sigma
^{X_t^{x,\varepsilon}} } $. Here, $\operatorname{adj} A$ is the adjugate matrix of
$A$. By the property of the Skorohod integral [Proposition 1.3.3 of
\citet{Nualart} and Lemma 5.2 of \citet{Malliavin} or (4.15) of proof of
Lemma 4.10 of \citet{Malliavin-Thalmaier}],
we have
\begin{eqnarray}\label{MalliavinIBP_Y}
\qquad&& H_{(i)}\bigl(X_t^{x,\varepsilon},G(t,x)\bigr)\nonumber
\\
&&\qquad =\delta
\Biggl(\sum_{j=1}^{N}G(t,x)
\gamma_{ij}^{X_t^{\varepsilon
}}DX_t^{x,\varepsilon,j} \Biggr)
\nonumber\\[-8pt]\\[-8pt]\nonumber
&&\qquad = \Biggl[G(t,x) \sum_{j=1}^{N} \sum
_{k=1}^d \int_0^t
\gamma_{ij}^{X_t^{\varepsilon}}\bigl(J_t^{x,\varepsilon
}\bigl(J_s^{x,\varepsilon}\bigr)^{-1}\varepsilon
V_k\bigl(X_s^{x,\varepsilon}\bigr)\bigr)^j
\,dB_s^k
\nonumber
\\
&&\quad\qquad{}\hspace*{4pt} -\sum_{j=1}^{N}\sum
_{k=1}^d \int_0^t
\bigl[D_{s,k} G(t,x)\bigr] \gamma_{ij}^{X_t^{\varepsilon
}}\bigl(J_t^{x,\varepsilon}\bigl(J_s^{x,\varepsilon}
\bigr)^{-1}\varepsilon V_k\bigl(X_s^{x,\varepsilon}
\bigr) \bigr)^j \,ds \Biggr].\nonumber
\end{eqnarray}
Again, by Lemma \ref{Properties}, the first and the second terms in
the second equality is characterized by
%
\begin{eqnarray}
G(t,x) \sum_{j=1}^{N} \sum
_{k=1}^d \int_0^t
\gamma_{ij}^{X_t^{\varepsilon}} \bigl(J_t^{x,\varepsilon
}
\bigl(J_s^{x,\varepsilon}\bigr)^{-1}\varepsilon
V_k\bigl(X_s^{x,\varepsilon}\bigr)\bigr)^j
\,dB_s^k &\in&{\mathcal K}_{r-1},
\\
\int_0^t \bigl[D_{s,k} G(t,x)
\bigr] \gamma_{ij}^{X_t^{\varepsilon
}}\bigl(J_t^{x,\varepsilon}
\bigl(J_s^{x,\varepsilon}\bigr)^{-1}\varepsilon
V_k\bigl(X_s^{x,\varepsilon}\bigr) \bigr)^j
\,ds &\in& {\mathcal K}_{r},
\end{eqnarray}
since $J_t^{x,\varepsilon}, (J_t^{x,\varepsilon})^{-1} \in{\mathcal
K}_0$, $\gamma_{ij}^{X_t^{\varepsilon}} \in{\mathcal K}_{-2}$ and
%
\begin{eqnarray}
&& \int_0^t \gamma_{ij}^{X_t^{\varepsilon}}
\bigl(J_t^{x,\varepsilon
}\bigl(J_s^{x,\varepsilon}
\bigr)^{-1}\varepsilon V_k\bigl(X_s^{x,\varepsilon}
\bigr)\bigr)^j \,dB_s^k \in{\mathcal
K}_{-2+1}={\mathcal K}_{-1}.
\end{eqnarray}
Then $H_{(i)} (X_t^{x,\varepsilon},G(t,x)) \in{\mathcal K}_{r-1}$ and
$H_{\alpha^{(k)}} (X_t^{x,\varepsilon},G(t,x)) \in{\mathcal
K}_{r-k}$. Therefore, we have the assertion.
\end{pf}

\section{Weak approximation with asymptotic expansion method}\label{sec3}
In the remainder of the paper, we use the following norms and seminorms:
%
\begin{eqnarray}
\llVert f \rrVert _{\infty} &=&\sup_{x \in{\mathbf R}^N} \bigl\llvert
f(x) \bigr\rrvert,\qquad \llVert \nabla f \rrVert _{\infty} = \max
_{i \in\{1,\ldots,N \}} \biggl\llVert \frac
{\partial f}{\partial x_i} \biggr\rrVert
_{\infty},
\\
\bigl\llVert \nabla^i f \bigr\rrVert _{\infty}&=& \max
_{j_1,\ldots,j_i \in\{1,\ldots,N
\}} \biggl\llVert \frac{\partial^i f}{\partial x_{j_1}\cdots\partial
x_{j_i} } \biggr\rrVert
_{\infty},\qquad f \in C_b^{\infty}\bigl(\R^N;\R
\bigr).
\end{eqnarray}
In the first step, we give approximation results of an asymptotic
expansion with Malliavin weights for $E[f(X_t^{x,\varepsilon})]$ where
%
\begin{eqnarray}
X_{t}^{x,\varepsilon}&=&x + \int_{0}^{t}V_{0}
\bigl(\varepsilon, X_{s}^{x,\varepsilon}\bigr) \,ds +\varepsilon\sum
_{j=1}^{d} \int_{0}^{t}V_{j}
\bigl(X_{s}^{x,\varepsilon}\bigr) \,dB_{s}^{j}.
\end{eqnarray}
Under the smoothness of the vector fields $V_j$, $j=0,1,\ldots,d$,
$X_{t}^{x,\varepsilon}$ is expanded as
%
\begin{eqnarray}
&& X_{t}^{x,\varepsilon}= X_{t}^{x,0}+\varepsilon
\frac{\partial}{\partial\varepsilon} X_{t}^{x,\varepsilon}\bigg| _{\varepsilon=0} +
\varepsilon^2 \frac{1}{2!} \frac{\partial^2}{\partial\varepsilon
^2}X_{t}^{x,\varepsilon}
\bigg| _{\varepsilon=0} +\cdots\qquad\mbox{in }{\mathbf D}^{\infty}.\label{Xexpansion}
\end{eqnarray}
Here, the above expansion in the space ${\mathbf D}^{\infty}$ is given in
the sense that for all $m \in{\mathbf N}$,
%
\begin{eqnarray}
&&\limsup_{\varepsilon\downarrow0} \frac{1}{\varepsilon^{m+1}} \Biggl\llVert
X_t^{x,\varepsilon}- \Biggl\{ X_{t}^{x,0}+\sum
_{i=1}^m \varepsilon^i
\frac{1}{i!}\frac{\partial^{i}}{\partial\varepsilon
^{i}}X_t^{x,\varepsilon}
\bigg| _{\varepsilon=0} \Biggr\} \Biggr\rrVert _{{\mathbf
D}^{k,p}} < \infty
\nonumber\\[-8pt]\\[-8pt]
\eqntext{\forall k \in{\mathbf N}, \forall p\in[1, \infty).}
\end{eqnarray}
For instance, see \citet{Watanabe} and \citet{KT} for
the details.

Let us define $\bar{X}_t^{x,\varepsilon}$ as the sum of the first two
terms in the expansion (\ref{Xexpansion}) as follows:
%
\begin{eqnarray}
\bar{X}_t^{x,\varepsilon}&=&X_{t}^{x,0}+
\varepsilon\frac{\partial
}{\partial\varepsilon} X_{t}^{x,\varepsilon}\bigg|
_{\varepsilon=0}.
\end{eqnarray}
We remark that $X_{t}^{x,0}$ is the solution to the following ODE:
%
\begin{eqnarray}
X_{t}^{x,0}&=&x + \int_{0}^{t}V_{0}
\bigl(0, X_{s}^{x,0}\bigr) \,ds,
\end{eqnarray}
and $\frac{\partial}{\partial\varepsilon}X_{s}^{x,\varepsilon
}\mid_{\varepsilon=0}$ satisfies the following linear SDE:
%
\begin{eqnarray}
\frac{\partial}{\partial\varepsilon}X_{s}^{x,\varepsilon,l}\bigg| _{\varepsilon=0} &=&\int
_{0}^{t} \frac{\partial}{\partial\varepsilon} V_{0}^l
\bigl(\varepsilon, X_{s}^{x,0}\bigr)\bigg| _{\varepsilon=0}\,ds+
\sum_{j=1}^{d}\int_{0}^{t}V_{j}^l
\bigl(X_{s}^{x,0}\bigr)\,dB_{s}^{j}
\nonumber\\[-8pt]\\[-8pt]\nonumber
&&{}+ \sum_{k=1}^N \int
_{0}^{t} \partial_k
V_{0}^l\bigl(0,X_{s}^{x,\varepsilon
}\bigr)
\bigg| _{\varepsilon=0}\frac{\partial}{\partial\varepsilon
}X_{s}^{x,\varepsilon,k}
\bigg| _{\varepsilon=0}\,ds,
\nonumber
\\
\frac{\partial}{\partial\varepsilon}X_{0}^{x,\varepsilon,l}\bigg|_{\varepsilon=0}&=&0,\qquad l=1,
\ldots,N.
\end{eqnarray}
The solution of $\frac{\partial}{\partial\varepsilon
}X_{s}^{x,\varepsilon}\mid _{\varepsilon=0}$ is given by
\begin{eqnarray}
&&\sum_{j=1}^d \int_0^t
J_t^{x,0} \bigl(J_u^{x,0}
\bigr)^{-1} V_j\bigl({X}_u^{x,0}
\bigr) \,dB_u^j
\nonumber\\[-8pt]\\[-8pt]\nonumber
&&\qquad{}  + \int_0^t
J_t^{x,0} \bigl(J_u^{x,0}
\bigr)^{-1} \frac{\partial}{\partial
\varepsilon} V_{0}\bigl(\varepsilon,
X_{u}^{x,0}\bigr)\bigg| _{\varepsilon
=0}\,du,
\nonumber
\end{eqnarray}
where $J_t^{x,0}=\nabla_x X_t^{x,0}$ [see (6.6) on page~354 of \citet{Kartzas-Shreve}, e.g.].
Note that $\frac{\partial}{\partial\varepsilon}
X_{t}^{x,\varepsilon}\mid _{\varepsilon=0}$ is a Gaussian random variable
with a mean $\mu(t)$ and a covariance matrix $\Sigma(t)=(\Sigma
_{i,j}(t))_{1\leq i,j \leq N}$
%
\begin{eqnarray}
\mu(t)&=&\int_0^t J_t^{x,0}
\bigl(J_u^{x,0}\bigr)^{-1} \frac{\partial
}{\partial\varepsilon}
V_{0}\bigl(\varepsilon, X_{u}^{x,0}\bigr)\bigg| _{\varepsilon =0}\,du,\label{mu}
\\
\Sigma_{i,j}(t)&=&\sum_{k=1}^d
\int_0^t \bigl(J_t^{x,0}
\bigl(J_s^{x,0}\bigr)^{-1} V_k
\bigl({X}_s^{x,0}\bigr)\bigr)^i
\bigl(J_t^{x,0} \bigl(J_s^{x,0}
\bigr)^{-1} V_k\bigl({X}_s^{x,0}\bigr)
\bigr)^j \,ds.\label{cov}\hspace*{-30pt}
\end{eqnarray}
Here, we note that $t \mapsto\mu(t)$ and $t \mapsto\Sigma
_{i,j}(t)$, $1\leq i,j \leq N$, are deterministic functions.
Therefore, $\bar{X}_t^{x,\varepsilon}
=X_{t}^{x,0}+\varepsilon\frac{\partial}{\partial\varepsilon}
X_{t}^{x,\varepsilon}\mid _{\varepsilon=0}$ is a Gaussian random variable
with a mean $X_{t}^{x,0}+\varepsilon\mu(t)$ and a covariance matrix
$\varepsilon^2 \Sigma(t)=(\varepsilon^2 \Sigma_{i,j}(t))_{1\leq i,j
\leq N}$.

\begin{rem}
1. When ${V}_0(\varepsilon,x)=\varepsilon{V}_0(x)$,
$\bar{X}_t^{x,\varepsilon}$ is given by
%
\begin{eqnarray}
\bar{X}_t^{x,\varepsilon} &=&x+\varepsilon\sum
_{i=0}^{d}V_{i}(x)\int
_{0}^{t}dB_{s}^{i},
\end{eqnarray}
where $B_t^0=t$.

2. When ${V}_0(\varepsilon,x)={V}_0(x)$,
$\bar{X}_t^{x,\varepsilon}$ is given by
%
\begin{eqnarray}
\bar{X}_t^{x,\varepsilon} &=&{X}_t^{x,0}+
\varepsilon\sum_{j=1}^d\int
_0^t J_t^{x,0}
\bigl(J_u^{x,0}\bigr)^{-1} V_j
\bigl({X}_u^{x,0}\bigr) \,dB_u^j.
\end{eqnarray}
\end{rem}

The next theorem shows the local approximation errors for
$E[f(X_t^{x,\varepsilon})]$ using Malliavin weights.

\begin{teo}\label{11111111111111}
Under condition \textup{[\textbf{H}]}, we have the following:
\begin{longlist}[2.]
\item[1.] For any $t \in(0,1]$ and $f \in C_b^{\infty}({\mathbf R}^N;{\mathbf R})$,
there exists $C>0$ such that
%
\begin{eqnarray}
&&\sup_{x \in{\mathbf R}^N} \Biggl\llvert E\bigl[f\bigl(X_t^{x,\varepsilon}
\bigr)\bigr]- \Biggl\{E \bigl[ f\bigl(\bar{X}_t^{x,\varepsilon}\bigr)
\bigr] +\sum_{j=1}^m
\varepsilon^jE\bigl[f\bigl(\bar {X}_t^{x,\varepsilon}
\bigr) \Phi_t^j\bigr] \Biggr\} \Biggr\rrvert
\nonumber\\[-8pt]\\[-8pt]\nonumber
&&\qquad\leq\varepsilon^{m+1} C \Biggl( \sum_{k=1}^{m+1}
t^{(m+1+k)/2}\bigl\llVert \nabla ^k f \bigr\rrVert _{\infty}
\Biggr),
\end{eqnarray}
where $\Phi_t^j$, $j \geq1$, is the Malliavin weights defined by
%
\begin{eqnarray}\label{Mw}
\Phi_t^j&=&\sum_{k=1}^j
\sum_{\beta_{1}+\cdots+\beta_{k}=j+k,
\beta_l \geq2} \sum_{\alpha^{(k)}\in\{1,\ldots,N \}^k}
\frac{1}{k!}
\nonumber\\[-8pt]\\[-8pt]\nonumber
&&{}\times H_{{\alpha}^{(k)}} \Biggl(\frac{\partial}{\partial
\varepsilon} X_{t}^{x,\varepsilon}
\bigg| _{\varepsilon=0}, \prod_{l=1}^k
\frac{1}{\beta_l !}\frac{\partial^{\beta_{l}}}{\partial
\varepsilon^{\beta_{l}}} X_t^{x,\varepsilon,\alpha
_{l}}\bigg|
_{\varepsilon=0} \Biggr).
\end{eqnarray}
\item[2.]
For any $t \in(0,1]$ and Lipschitz continuous function $f\dvtx  \R^N
\rightarrow\R$, there exists $C>0$ such that
%
\begin{eqnarray}
&& \sup_{x \in{\mathbf R}^N} \Biggl\llvert E\bigl[f\bigl(X_t^{x,\varepsilon}
\bigr)\bigr]- \Biggl\{E \bigl[ f\bigl(\bar{X}_t^{x,\varepsilon}\bigr)
\bigr] +\sum_{j=1}^m
\varepsilon^jE\bigl[f\bigl(\bar {X}_t^{x,\varepsilon}
\bigr) \Phi_t^j\bigr] \Biggr\} \Biggr\rrvert
\nonumber\\[-8pt]\\[-8pt]\nonumber
&&\qquad \leq
\varepsilon ^{m+1} C t^{(m+2)/2},
\end{eqnarray}
with same weights in $($\ref{Mw}$)$.
\item[3.]
For any $t \in(0,1]$ and bounded Borel function $f\dvtx  \R^N \rightarrow
\R$, there exists \mbox{$C>0$} such that
%
\begin{eqnarray}
&& \sup_{x \in{\mathbf R}^N} \Biggl\llvert E\bigl[f\bigl(X_t^{x,\varepsilon}
\bigr)\bigr]- \Biggl\{E \bigl[ f\bigl(\bar{X}_t^{x,\varepsilon}\bigr)
\bigr] +\sum_{j=1}^m
\varepsilon^jE\bigl[f\bigl(\bar {X}_t^{x,\varepsilon}
\bigr) \Phi_t^j\bigr] \Biggr\} \Biggr\rrvert
\nonumber\\[-8pt]\\[-8pt]\nonumber
&&\qquad \leq \varepsilon ^{m+1} C t^{(m+1)/2},
\end{eqnarray}
with same weights in $($\ref{Mw}$)$.
\end{longlist}
\end{teo}

\begin{pf}
See Appendix~\ref{appeA}.
\end{pf}

\begin{rem}
When $\tilde{V}_0(\varepsilon,x)=\varepsilon\tilde{V}_0(x)$,
$X_t^{x,\varepsilon}$ has the following expansion:
\begin{eqnarray*}
X_t^{x,\varepsilon} &=&x+\varepsilon\sum
_{j=0}^{d}\tilde{V}_{j}(x)\int
_{0}^{t} \circ dB_{s}^{j}
\\
&&{}+\sum_{k=2}^m \varepsilon^k
\sum_{(i_1,\ldots,i_k)\in\{
0,1,\ldots,d \}^k} (\tilde{V}_{i_1}\cdots\tilde {V}_{i_k}) (x)
\int
_{0<t_1<\cdots<t_k<t}\circ dB_{t_1}^{i_1}\circ \cdots\circ
dB_{t_k}^{i_k}
\nonumber
\\
&&{}+\varepsilon^{m+1} \tilde{R}_m(t,x,\varepsilon),
\end{eqnarray*}
where $\tilde{R}_m(t,x,\varepsilon)$ is the residual. Here, we used
the notation $B_t^0=t$. Then
\begin{eqnarray*}
&& \frac{1}{k!}\frac{\partial^k}{\partial\varepsilon^k} X_t^{x,\varepsilon}\bigg|
_{\varepsilon=0}
\\
&&\qquad =\sum_{(i_1,\ldots,i_k)\in\{0,1,\ldots,d \}^k}
(\tilde{V}_{i_1}\cdots\tilde{V}_{i_k}) (x) \int
_{0<t_1<\cdots<t_k<t}\circ dB_{t_1}^{i_1}\circ\cdots\circ
dB_{t_k}^{i_k}.
\end{eqnarray*}
\end{rem}

\begin{rem}
$\Phi_t^j$ is obtained by multiple Skorohod integral
and each Malliavin weight is concretely calculated as follows; for
$G(t,x) \in{\mathcal K}_r$ and $i=1,\ldots,N$,
%
\begin{eqnarray}
\qquad &&H_{(i)} \biggl(\frac{\partial}{\partial\varepsilon} X_{t}^{x,\varepsilon}
\bigg| _{\varepsilon=0}, G(t,x) \biggr)
\nonumber
\\
&&\qquad =G(t,x) \sum_{j=1}^{N}\sum
_{k=1}^{d} \bigl[\Sigma(t)^{-1}
\bigr]_{i,j} \int_0^t
\bigl(J_t^{x,0} \bigl(J_s^{x,0}
\bigr)^{-1} V_k\bigl({X}_s^{x,0}\bigr)
\bigr)^j \,dB_s^k
\\
&&\quad\qquad{}-\sum_{j=1}^{N}\sum
_{k=1}^{d} \bigl[\Sigma(t)^{-1}
\bigr]_{i,j} \int_0^t
D_{s,k} G(t,x) \bigl(J_t^{x,0}
\bigl(J_s^{x,0}\bigr)^{-1} V_k
\bigl({X}_s^{x,0}\bigr)\bigr)^j \,ds
\nonumber
\end{eqnarray}
with the deterministic covariance matrix $(\Sigma_{i,j}(t))_{1\leq i,j
\leq N}$ corresponds to $($\ref{cov}$)$, {that is},
\begin{eqnarray}
\qquad\Sigma_{i,j}(t)&=&\sum_{k=1}^d
\int_0^t D_{s,k}\frac{\partial
}{\partial\varepsilon}
X_{t}^{x,\varepsilon,i}\bigg| _{\varepsilon
=0}D_{s,k}
\frac{\partial}{\partial\varepsilon} X_{t}^{x,\varepsilon,j}\bigg| _{\varepsilon=0}\,ds
\nonumber\\[-8pt]\\[-8pt]\nonumber
&=&\sum_{k=1}^d \int
_0^t \bigl(J_t^{x,0}
\bigl(J_s^{x,0}\bigr)^{-1} V_k
\bigl({X}_s^{x,0}\bigr)\bigr)^i
\bigl(J_t^{x,0} \bigl(J_s^{x,0}
\bigr)^{-1} V_k\bigl({X}_s^{x,0}\bigr)
\bigr)^j \,ds.
\nonumber
\end{eqnarray}
\end{rem}

Let $(P_t)_t$ be linear operators on $f \in C_b(\R^N;\R)$ defined by
%
\begin{equation}
P_t f(x)=E\bigl[f\bigl(X_t^{x,\varepsilon}\bigr)\bigr].
\end{equation}
We remark that $(P_t)_t$ is a semigroup. Also let $(\bar{P}_t)_t$ be
linear operators on $f \in C_b(\R^N;\R)$ defined by
%
\begin{equation}
\bar{P}_t f(x)=E\bigl[f\bigl(\bar{X}_t^{x,\varepsilon}
\bigr)\bigr].
\end{equation}

Next, as an approximation of $P_s$ we introduce a linear operator
$Q^m_{(s)}$ below.
First, for $j\geq1$ and $t \in(0,1]$,
let $\bar{P}_{\Phi^j}(t)$ be a linear operator defined by the
following expectation with Malliavin weight $\Phi_t^j$:
%
\begin{equation}
\bar{P}_{\Phi^j}(t)f(x)=E \bigl[f \bigl( \bar{X}_t^{x,\varepsilon}
\bigr)\Phi _t^j \bigr].
\end{equation}
Then $( Q^m_{(s)})_{s \in(0,1]}$ is defined as linear operators:
%
\begin{eqnarray}
Q^m_{(s)}f(x)&=&\bar{P}_s f(x)+\sum
_{j=1}^{m} \varepsilon^{j} \bar
{P}_{\Phi^j}(s)f(x).
\end{eqnarray}

We remark that
%
\begin{eqnarray}
\bar{P}_{\Phi^j}(t)f(x) &=&\int_{\R^N}f(y)E \bigl[
\Phi_t^j \mid \bar{X}_t^{x,\varepsilon}=y
\bigr] p^{\bar{X}^\varepsilon}(t,x,y)\,dy
\\
&=&E \bigl[f \bigl( \bar{X}_t^{x,\varepsilon}\bigr) {\mathcal
M}_{(j)}\bigl(t,x,\bar {X}_t^{x,\varepsilon}\bigr) \bigr],
\end{eqnarray}
where ${\mathcal M}_{(j)}(t,x,y)=E [ \Phi_t^j \mid  \bar
{X}_t^{x,\varepsilon}=y ]$ and $y \mapsto p^{\bar{X}^\varepsilon
}(t,x,y)$ is the density of~$\bar{X}_t^{x,\varepsilon}$.

Then $Q^m_{(s)}$ can be written as follows:
%
\begin{eqnarray}
Q^m_{(s)}f(x)&=&E\bigl[f\bigl(\bar{X}_s^{x,\varepsilon}
\bigr){\mathcal M}^{m}\bigl(s,x,\bar{X}_s^{x,\varepsilon}
\bigr) \bigr],
\end{eqnarray}
where
${\mathcal M}^{m}(s,x,y)=1+\sum_{j=1}^m \varepsilon^{j} {\mathcal
M}_{(j)}(s,x,y)$.

Then we have the following explicit representation for the Malliavin
weight function ${\mathcal M}^m$.

\begin{teo}\label{222222222222222222}\label{Mallfunction}
Under condition \textup{[\textbf{H}]}, the Malliavin weight function ${\mathcal
M}^m$ is given by
%
\begin{eqnarray}\label{formulaMweight}
&&{\mathcal M}^{m}(t,x,y)\nonumber
\\[-1pt]
&&\qquad =1+\sum_{j=1}^m \varepsilon^{j}
\sum_{k=1}^j \sum
_{\beta_{1}+\cdots+\beta_{k}=j+k, \beta_l \geq2} \sum_{\alpha^{(k)}\in\{1,\ldots,N \}^k}
\frac{\varepsilon^k}{k!}
\\[-1pt]
&&\hspace*{49pt}{}\times \partial^{\ast}_{\alpha_k} \circ\partial^{\ast}_{\alpha
_{k-1}}
\circ\cdots\circ\partial^{\ast}_{\alpha_1} E \Biggl[ \prod
_{l=1}^k \frac{1}{\beta_l !}\frac{\partial^{\beta
_{l}}}{\partial\varepsilon^{\beta_{l}}}
X_t^{x,\varepsilon,\alpha
_{l}}\bigg| _{\varepsilon=0} \Big| \bar{X}_{t}^{x,\varepsilon}=y
\Biggr],\nonumber\hspace*{-15pt}
\end{eqnarray}
where $\partial^{\ast}$ is the divergence operator on the Gaussian
space $({\mathbf R}^N,\nu)$, {that is},
%
\begin{eqnarray}
\nu(dy)&=&p^{\bar{X}^\varepsilon}(t,x,y)\,dy
\nonumber
\\[-1pt]
&=& \frac{1}{ (2\pi\varepsilon)^{N/2} \det(\Sigma(t))^{1/2}}\nonumber
\\[-1pt]
&&{}\times  e^{-\afrac{(y-X_{t}^{x,0}-\varepsilon\mu(t))^{\top} \Sigma
^{-1}(t)(y-X_{t}^{x,0}-\varepsilon\mu(t)) }{2\varepsilon
^2}}\,dy,
\\[-1pt]
\partial^{\ast}_i A(y)&=&- \biggl[\frac{\partial}{\partial y_i}\log
p^{\bar{X}^\varepsilon}(t,x,y) \biggr] A(y)- \frac{\partial
}{\partial y_i}A(y),\nonumber
\\[-1pt]
\eqntext{ A \in{\mathcal S}
\bigl({\mathbf R}^N\bigr), 1\leq i \leq N.}
\end{eqnarray}
Here, $\mu(t)$ and $\Sigma(t)=(\Sigma_{i,j}(t))_{1\leq i,j \leq N}$ are
defined in $($\ref{mu}$)$ and $($\ref{cov}$)$, respectively, that is,
%
\begin{eqnarray}
\mu(t)&=&\int_0^t J_t^{x,0}
\bigl(J_u^{x,0}\bigr)^{-1} \frac{\partial
}{\partial\varepsilon}
V_{0}\bigl(\varepsilon, X_{u}^{x,0}\bigr)\bigg|
_{\varepsilon
=0}\,du,
\\
\Sigma_{i,j}(t) &=&\sum_{k=1}^d
\int_0^t \bigl(J_t^{x,0}
\bigl(J_s^{x,0}\bigr)^{-1} V_k
\bigl({X}_s^{x,0}\bigr)\bigr)^i
\bigl(J_t^{x,0} \bigl(J_s^{x,0}
\bigr)^{-1} V_k\bigl({X}_s^{x,0}\bigr)
\bigr)^j \,ds,\hspace*{-40pt}
\end{eqnarray}
and ${\mathcal S}({\mathbf R}^N)$ is the Schwartz rapidly decreasing
functions on ${\mathbf R}^N$.
\end{teo}

\begin{pf}
See Appendix~\ref{appeB}.
\end{pf}

\begin{rem}
The\vspace*{1pt} term
$\prod_{l=1}^k \frac{1}{\beta_l !}\frac{\partial^{\beta
_{l}}}{\partial\varepsilon^{\beta_{l}}} X_t^{x,\varepsilon,\alpha
_{l}}\mid _{\varepsilon=0}$
in each conditional expectation in (\ref{formulaMweight}) of Theorem~\ref{Mallfunction}
is generally expressed as a finite sum of iterated multiple Wiener--It\^o integrals.
Hence, we are able to explicitly compute each conditional expectation,
conditioned on $\bar{X}_{t}^{x,\varepsilon}$ that is given
by the first-order Wiener--It\^o
integral.

For instance, let
$q_k =({q}_{k,1},\ldots,{q}_{k,d})^{\top}$, ${q}_{k,i} \in L^2
([0,t])$, $k=1,2,3,4$, $i=1,\ldots,d$ and
$h_l(\xi;v)$ be the (one dimensional) Hermite polynomial of degree $l$
with parameter
$v= \int_0^t q^{\top}_{1}(s)q_{1}(s) \,ds$.
Then the conditional expectations of the second- and the third-order
iterated multiple Wiener--It\^o
integrals are evaluated
as the following formulas:
\begin{eqnarray}
&& \E \biggl[\int_0^t \int
_0^s q^{\top}_{2}(u)\,dB_uq^{\top}_{3}(s)\,dB_s
\Big| \int_0^t q^{\top}_{1}(s)\,dB_s=
\xi \biggr]
\nonumber\\[-8pt]\\[-8pt]\nonumber
&&\qquad = \biggl(\int_0^t \int
_0^s q^{\top}_{2}(u)q_{1}(u)\,du\,
q^{\top
}_{3}(s)q_{1}(s)\,ds \biggr)
\frac{h_2(\xi;v)}{v^2},
\nonumber
\\
&& \E \biggl[\int_0^t \int
_0^s\int_0^u
q^{\top}_{2}(r)\,dB_r q^{\top}_{3}(u)\,dB_u
q^{\top}_{4}(s)\,dB_s \Big| \int
_0^t q^{\top}_{1}(s)\,dB_s=
\xi \biggr]
\nonumber\\[-8pt]\\[-8pt]\nonumber
&&\qquad = \biggl(\int_0^t q^{\top}_{4}(s)q_{1}(s)
\int_0^s q^{\top}_{3}(u)q_{1}(u)
\int_0^u q^{\top}_{2}(r)q_{1}(r)\,dr
\,du\,ds \biggr) \frac{h_3(\xi;v)}{v^3},\nonumber\hspace*{-25pt}
\end{eqnarray}
where
$h_2(\xi;v) = \xi^2 -v$ and $h_3(\xi;v) = \xi^3 -3 v\xi$.

The conditional expectations of higher order iterated multiple
Wiener--It\^o integrals
can be evaluated in the similar manner.
For the details, see \citet{Takahashi} and \citet{Takahashi-Takehara-Todaa}.
In fact, we obtain the Malliavin weights appearing in the numerical
examples in Section~\ref{sec5}
as closed forms by applying the formulas.
\end{rem}

Therefore, Theorem~\ref{11111111111111} is summarized as follows.

\begin{cor}\label{shorttimeerror}
Assume that condition \textup{[\textbf{H}]} holds.
\begin{longlist}[2.]
\item[1.]
There exists $C>0$ such that
%
\begin{equation}
\label{shorttimeerror1} \bigl\llVert P_sf-Q^m_{(s)}f
\bigr\rrVert _{\infty}\leq\varepsilon^{m+1} C \Biggl( \sum
_{k=1}^{m+1} s^{(m+1+k)/2}\bigl\llVert
\nabla^k f \bigr\rrVert _{\infty} \Biggr),
\end{equation}
for any $s \in(0,1]$ and $f \in C_b^{\infty}({\mathbf R}^N;{\mathbf R})$.
\item[2.]
There exists $C>0$ such that
%
\begin{equation}
\label{shorttimeerror2} \bigl\llVert P_sf-Q^m_{(s)}f
\bigr\rrVert _{\infty}\leq\varepsilon^{m+1} C s^{(m+2)/2},
\end{equation}
for any $s \in(0,1]$ and Lipschitz continuous function $f\dvtx  \R^N
\rightarrow\R$.
\item[3.]
There exists $C>0$ such that
%
\begin{equation}
\label{shorttimeerror3} \bigl\llVert P_sf-Q^m_{(s)}f
\bigr\rrVert _{\infty}\leq\varepsilon^{m+1} C s^{(m+1)/2},
\end{equation}
for any $s \in(0,1]$ and bounded Borel function $f\dvtx  \R^N \rightarrow
\R$.
\end{longlist}
\end{cor}
%

\begin{rem}
The above results are obtained based on the integration by parts
argument for $G(s,x)\in{\mathcal K}_r$ with time $s \in(0,1]$.
However, we are able to show
that the same results hold for $s \in(0,T]$, $T>0$, using the
properties of the elements in the space ${\mathcal K}_r^T$ defined as
in \citet{Crisan-Manolarakis-Nee}.
\end{rem}

Next,
for $T>0$, $\gamma>0$, define a partition $\pi=\{ (t_0,t_1,\ldots,t_n)\dvtx  0=t_0<t_1<\cdots<t_n=T,t_k=k^{\gamma}T/n^{\gamma}, n\in{\mathbf
N} \}$ and $s_k=t_k-t_{k-1}$, $k=1,\ldots,n$. Using the asymptotic
expansion operator $Q^m$ of $P$, we can guess the following semigroup
approximation.
\begin{eqnarray*}
E\bigl[ f\bigl(X_T^{x,\varepsilon}\bigr) \bigr]=P_Tf(x)=P_{s_n}P_{s_{n-1}}
\cdots P_{s_1} f(x) &\simeq& Q^m_{(s_n)}
Q^m_{(s_{n-1})}\cdots Q^m_{(s_1)} f(x).
\end{eqnarray*}

The next theorem shows our main result on the approximation error for
this scheme.

\begin{teo}\label{33333333333333}
Assume that condition \textup{[\textbf{H}]} holds.
Let $T>0$, $\gamma>0$ and $n \in\N$.
\begin{longlist}[3.]
\item[1.]
For any $f \in C_b^{\infty}({\mathbf R}^N;{\mathbf R})$, there exists $C>0$
such that
%
\begin{eqnarray}
&&\bigl\llVert P_Tf-Q^m_{(s_n)}
Q^m_{(s_{n-1})}\cdots Q^m_{(s_1)} f \bigr
\rrVert _{\infty
}\leq\varepsilon^{m+1} \frac{C}{n^{\gamma(m+2)/2}},
\nonumber\\[-8.5pt]\\[-8.5pt]
\eqntext{0<\gamma < m/(m+2),}
\\
&&\bigl\llVert P_Tf-Q^m_{(s_n)}
Q^m_{(s_{n-1})}\cdots Q^m_{(s_1)} f \bigr
\rrVert _{\infty
}\leq\varepsilon^{m+1} \frac{C}{n^{m/2}}(1+\log
n),
\nonumber\\[-8.5pt]\\[-8.5pt]
\eqntext{\gamma =m/(m+2),}
\\
&&\bigl\llVert P_Tf-Q^m_{(s_n)}
Q^m_{(s_{n-1})}\cdots Q^m_{(s_1)} f \bigr
\rrVert _{\infty
}\leq\varepsilon^{m+1} \frac{C}{n^{m/2}},
\qquad
\gamma>m/(m+2).\hspace*{-30pt}
\end{eqnarray}
\item[2.]
For any Lipschitz continuous function $f\dvtx  \R^N \rightarrow\R$,
there exists $C>0$ such that
%
\begin{eqnarray}
&&\bigl\llVert P_Tf-Q^m_{(s_n)}
Q^m_{(s_{n-1})}\cdots Q^m_{(s_1)} f \bigr
\rrVert _{\infty
}\leq\varepsilon^{m+1} \frac{C}{n^{\gamma(m+2)/2}},
\nonumber\\[-8.5pt]\\[-8.5pt]
\eqntext{0<\gamma < m/(m+2),}
\\
&&\bigl\llVert P_Tf-Q^m_{(s_n)}
Q^m_{(s_{n-1})}\cdots Q^m_{(s_1)} f \bigr
\rrVert _{\infty
}\leq\varepsilon^{m+1} \frac{C}{n^{m/2}}(1+\log
n),
\nonumber\\[-8.5pt]\\[-8.5pt]
\eqntext{\gamma =m/(m+2),}
\\
&&\bigl\llVert P_Tf-Q^m_{(s_n)}
Q^m_{(s_{n-1})}\cdots Q^m_{(s_1)} f \bigr
\rrVert _{\infty
}\leq\varepsilon^{m+1} \frac{C}{n^{m/2}},
\qquad \gamma> m/(m+2).\hspace*{-30pt}
\end{eqnarray}
\item[3.]
For any bounded Borel function $f\dvtx  \R^N \rightarrow\R$, there
exists $C>0$ such that
%
\begin{eqnarray}
&&\bigl\llVert P_Tf-Q^m_{(s_n)}
Q^m_{(s_{n-1})}\cdots Q^m_{(s_1)} f \bigr
\rrVert _{\infty
}\leq\varepsilon^{m+1} \frac{C}{n^{\gamma(m+1)/2}},
\nonumber\\[-8.5pt]\\[-8.5pt]
\eqntext{0<\gamma< (m-1)/(m+1),}
\\
&&\bigl\llVert P_Tf-Q^m_{(s_n)}
Q^m_{(s_{n-1})}\cdots Q^m_{(s_1)} f \bigr
\rrVert _{\infty
}\leq\varepsilon^{m+1} \frac{C}{n^{(m-1)/2}}(1+\log
n),
\nonumber\\[-8.5pt]\\[-8.5pt]
\eqntext{\gamma =(m-1)/(m+1),}
\\
&&\bigl\llVert P_Tf-Q^m_{(s_n)}
Q^m_{(s_{n-1})}\cdots Q^m_{(s_1)} f \bigr
\rrVert _{\infty
}\leq\varepsilon^{m+1} \frac{C}{n^{(m-1)/2}},
\nonumber\\[-8.5pt]\\[-8.5pt]
\eqntext{\gamma>(m-1)/(m+1).}
\end{eqnarray}
\end{longlist}
\end{teo}

\begin{pf}
See Appendix~\ref{appeC}.
\end{pf}

\begin{rem}
Due to the theorem above, the higher order asymptotic expansion
provides the higher order weak approximation.
In fact, we can mostly attain enough accuracy
even when the expansion order $m$ is low such as $m=1,2$.
In Section~\ref{sec5}, we confirm this fact through numerical examples.
\end{rem}

\begin{rem}
When $\gamma=1$, {that is}, $s_k=T/n$ for all $k=1,\ldots,n$, we have:
\begin{longlist}[3.]
\item[1.]
For any $f \in C_b^{\infty}({\mathbf R}^N;{\mathbf R})$, there exists $C>0$
such that
\begin{eqnarray}
&&\bigl\llVert P_Tf- \bigl(Q^m_{(T/n)}
\bigr)^n f \bigr\rrVert _{\infty}\leq\varepsilon^{m+1}
\frac
{C}{n^{m/2}}.
\nonumber
\end{eqnarray}
\item[2.]
For any Lipschitz continuous function $f\dvtx  \R^N \rightarrow\R$,
there exists $C>0$ such that
\begin{eqnarray}
&&\bigl\llVert P_Tf-\bigl(Q^m_{(T/n)}
\bigr)^n f \bigr\rrVert _{\infty}\leq\varepsilon^{m+1}
\frac
{C}{n^{m/2}}.
\nonumber
\end{eqnarray}
\item[3.]
For any bounded Borel function $f\dvtx  \R^N \rightarrow\R$, there
exists $C>0$ such that
\begin{eqnarray}
&&\bigl\llVert P_Tf-\bigl(Q^m_{(T/n)}
\bigr)^n f \bigr\rrVert _{\infty}\leq\varepsilon^{m+1}
\frac
{C}{n^{(m-1)/2}}.
\nonumber
\end{eqnarray}
\end{longlist}
\end{rem}
%

\section{Computation with Malliavin weights}\label{sec4}
This section illustrates computational scheme for implementation of our method.
\subsection{Backward discrete-time approximation}\label{sec4.1}
For preparation, we describe a backward discrete-time approximation of
our method.

For $s \in(0,1]$ and $x,y\in\R^N$, define $p^{m}(s,x,y)$ as
%
\begin{equation}
Q^m_{(s)}f(x)=\int_{\R^N}
f(y)p^{m}(s,x,y)\,dy.
\end{equation}
Then $p^{m}(s,x,y)$ is given by using the Malliavin weight function
${\mathcal M}^m$ as follows:
%
\begin{equation}
p^{m}(s,x,y)= {\mathcal M}^m(s,x,y)p^{\bar{X}^\varepsilon}(s,x,y),
\end{equation}
with
%
\begin{eqnarray}
p^{\bar{X}^\varepsilon}(s,x,y)&=& \frac{1}{ (2\pi\varepsilon
^2)^{N/2} \det(\Sigma(s))^{1/2}}
\nonumber\\[-8pt]\\[-8pt]\nonumber
&&{} \times e^{-\afrac{(y-\varepsilon\mu(s)-{X}_s^{x,0})^{\top} \Sigma
^{-1}(s)(y-\varepsilon\mu(s)-{X}_s^{x,0}) }{2 \varepsilon
^2}},
\end{eqnarray}
where $\mu(s)$ and $\Sigma(s)=(\Sigma_{i,j}(s))_{1\leq i,j \leq N}$
are defined in $($\ref{mu}$)$ and $($\ref{cov}$)$, respectively.

Then we are able to calculate $(Q^m_{(T/n)})^n f(x)$ as follows:
%
\begin{eqnarray}
&&\bigl(Q^m_{(T/n)}\bigr)^n f(x)
\\
&&\qquad = \int_{(\R^N)^n} f(y_n)\prod
_{i=0}^{n-1} p^m (s_{i},y_{i},y_{i+1})\,dy_{n}
\cdots dy_{1}
\\
&&\qquad =\int_{(\R^N)^{n-1}} q_{n-1}(y_{n-1}) \prod
_{i=0}^{n-2} p^m
(s_{i},y_{i},y_{i+1})\,dy_{n-1} \cdots
dy_{1}
\\
&&\qquad =\int_{(\R^N)^{n-2}} q_{n-2}(y_{n-2}) \prod
_{i=0}^{n-3} p^m
(s_{i},y_{i},y_{i+1})\,dy_{n-2} \cdots
dy_{1}
\\
&&\qquad =\int_{\R^N} q_{1}(y_{1})
p^m (s_{1},y_{0},y_{1})
\,dy_{1},
\end{eqnarray}
with $y_0=x$.

\subsection{Example of computational scheme}\label{sec4.2}
We are able to compute the expectation in the various ways such as
numerical integration
and Monte Carlo simulation.
As an illustrative purpose and an example,
this subsection briefly describes a scheme based on Monte Carlo simulation.

In computation of $(Q^m_{(T/n)})^n f(x)$ with simulation (for the case
of $\gamma=1$),
we store $\bar{X}^{x,(j)}_{T/n}\equiv\bar{X}^{x,\varepsilon,(j)}_{T/n}$, which stands for the $j$th ($1\leq j \leq M$) independent
outcome of $\bar{X}^{x,\varepsilon}$ at $T/n$
(i.e., at $t_i+T/n$)
starting from $x$ at each grid $t_i=(iT)/n$ ($0 \leq i \leq n-1$).

Then we calculate an approximate semigroup at each time grid.
That is, $q_{n-1}(x)$, $q_{n-2}(x)$ are calculated as follows:
%
\begin{eqnarray}
q_{n-1}(x) &=&\int_{\R^N}f(y)p^m(T/n,x,y)\,dy
\\
&=&\int_{\R^N}f(y){\mathcal M}^m(T/n,x,y)p^{\bar{X}^\varepsilon
}(T/n,x,y)\,dy
\\
&\simeq&\frac{1}{M}\sum_{j=1}^{M}
f\bigl(\bar{X}^{x,(j)}_{T/n}\bigr){\mathcal M}^m
\bigl(T/n,x,\bar{X}^{x,(j)}_{T/n}\bigr),
\\
q_{n-2}(x) &=&\int_{\R^N}q_{n-1}(y)p^m(T/n,x,y)\,dy
\\
&=&\int_{\R^N}q_{n-1}(y){\mathcal
M}^m(T/n,x,y) p^{\bar
{X}^\varepsilon}(T/n,x,y)\,dy
\\
&\simeq&\frac{1}{M}\sum_{j=1}^{M}
q_{n-1}\bigl(\bar {X}^{x,(j)}_{T/n}\bigr){\mathcal
M}^m\bigl(T/n,x,\bar{X}^{x,(j)}_{T/n}\bigr).
\end{eqnarray}
Therefore, in general,
%
\begin{eqnarray}
q_{i-1}(x) &=&\int_{\R^N}q_{i}(y)p^m(T/n,x,y)\,dy
\\
&=&\int_{\R^N}q_{i}(y){\mathcal
M}^m(T/n,x,y) p^{\bar{X}^\varepsilon
}(T/n,x,y)\,dy \label{ni1}
\\
&\simeq&\frac{1}{M}\sum_{j=1}^{M}
q_{i}\bigl(\bar {X}^{x,(j)}_{T/n}\bigr){\mathcal
M}^m\bigl(T/n,x,\bar{X}^{x,(j)}_{T/n}\bigr).
\label{MC1}
\end{eqnarray}
Finally, we obtain an approximation:
%
\begin{eqnarray}
\bigl(Q^m_{(T/n)}\bigr)^n f(x) &=&\int
_{\R^N}q_{1}(y)p^m(T/n,x,y)\,dy
\\
&=&\int_{\R^N}q_{1}(y){\mathcal
M}^m(T/n,x,y) p^{\bar{X}^\varepsilon
}(T/n,x,y)\,dy \label{ni2}
\\
&\simeq&\frac{1}{M}\sum_{j=1}^{M}
q_{1}\bigl(\bar {X}^{x,(j)}_{T/n}\bigr){\mathcal
M}^m\bigl(T/n,x,\bar{X}^{x,(j)}_{T/n}\bigr).
\label{MC2}
\end{eqnarray}
We also remark that if the numerical integration method is applied,
the scheme is based on equations (\ref{ni1}) and (\ref{ni2}).

\subsection{Comparison with Kusuoka--Lyons--Victoir (KLV) cubature method}\label{sec4.3}
In this subsection, we compare our method to a related work,
Kusuoka--Lyons--Victoir (KLV)
cubature method on Wiener space
[\citeauthor{Kusuokaa} (\citeyear{Kusuokaa,Kusuokad}), \citet{Lyons-Victoir}].

As mentioned\vspace*{1pt} above, we defined the operator $Q_{(s)}^m$
by using the asymptotic expansion with Malliavin weights, while \citeauthor{Kusuokaa}
(\citeyear{Kusuokaa,Kusuokad}) and \citet{Lyons-Victoir} developed a construction
method of a local approximation operator $\hat{Q}_{(s)}^{m}$ for $P_s$
based on finite variation paths $\omega_1,\ldots,\omega_l$ for some
$l \in{\mathbf N}$ with weights $\lambda_1,\ldots,\lambda_l$.

In the following, we summarize
our weak approximation method and the KLV cubature scheme.

\subsubsection*{Weak approximation with asympotic expansion and Malliavin weights}

Let $X_t^{x,\varepsilon}$ be a solution to the following SDE:
%
\begin{eqnarray}
 dX_t^{x,\varepsilon} &=& V_0\bigl(
\varepsilon,X_t^{x,\varepsilon
}\bigr)\,dt+ \varepsilon\sum
_{i=1}^d V_i\bigl(X_t^{x,\varepsilon}
\bigr) \,dB_t^i,\qquad X_0^{x,\varepsilon} = x.
\end{eqnarray}
For a Lipschitz continuous function $f$, $P_t
f(x)=E[f(X_t^{x,\varepsilon})]$ is approximated by
\vspace*{1pt}$Q_{(t)}^{m}f(x)=E[f(\bar{X}_t^{x,\varepsilon} )]+\sum_{j=1}^m
\varepsilon^j E[f(\bar{X}_t^{x,\varepsilon} ) \Phi_t^j ]=E[f(\bar
{X}_t^{x,\varepsilon} ){\mathcal M}^m(t,x,
\bar{X}_t^{x,\varepsilon})]$ as follows:
%
\begin{equation}
\bigl\llVert P_t f - Q_{(t)}^{m}f \bigr\rrVert
_{\infty}=O\bigl( \varepsilon^{m+1} t^{(m+2)/2}\bigr),\qquad t
\in(0,1].
\end{equation}
Then we have the global approximation,
%
\begin{equation}
\bigl\llVert P_T f - \bigl(Q_{(T/n)}^{m}
\bigr)^n f \bigr\rrVert _{\infty}=O\bigl( \varepsilon^{m+1}
n^{-m/2}\bigr).
\end{equation}

It is emphasized that we are able to evaluate Malliavin weights\break
${\mathcal M}^m(t,x,\bar{X}_t^{x,\varepsilon})$
mostly as closed forms
by applying computational schemes such as
conditional expectation formulas in \citet{Takahashi} and \citet{Takahashi-Takehara-Todaa}.
In fact, this is the case for the numerical examples in Section~\ref{sec5} of
this paper.

\subsubsection*{KLV cubature scheme on Wiener space}

Let $X_t^{x}$ be a solution to the following SDE:
%
\begin{eqnarray}
&&dX_t^x = V_0\bigl(X_t^x
\bigr)\,dt+\sum_{i=1}^d V_i
\bigl(X_t^x\bigr) \circ dB_t^i,\qquad
X_0^x = x.
\end{eqnarray}
A set of finite variation paths
$\omega=(\omega_1,\ldots,\omega_l)$ with $\lambda=(\lambda
_1,\ldots,\lambda_l)$ forms
\textit{cubature formula on Wiener space of degree} $m$ if for any $\alpha
\in{\mathcal A}_m$,
%
\begin{eqnarray}
&& E \biggl[\int_{0<t_1<\cdots<t_r<t} \circ dB_{t_1}^{\alpha_1}
\circ\cdots\circ dB_{t_r}^{\alpha
_r} \biggr]
\nonumber\\[-8pt]\\[-8pt]\nonumber
&&\qquad =\sum
_{j=1}^l \lambda_j \int
_{0<t_1<\cdots
<t_r<t}d\omega_{j,t_1}^{\alpha_1} \cdots d
\omega_{j,t_r}^{\alpha_r}.
\end{eqnarray}
$\omega=(\omega_1,\ldots,\omega_l)$ and $\lambda=(\lambda
_1,\ldots,\lambda_l)$ are called the cubature paths and weights, respectively.
Here, ${\mathcal A}_m$ is a set defined by ${\mathcal A}_m =\{ (\alpha
_1,\ldots,\alpha_r) \in\{ 0,1,\ldots, d \}^r;\break  r+\#\{j \mid  \alpha
_j=0 \}\leq m,r \in{\mathbf N} \}$.
For cubature paths $\omega=(\omega_1,\ldots,\omega_l)$ and weights
$\lambda=(\lambda_1,\ldots,\lambda_l)$, consider the following ODEs:
%
\begin{eqnarray}
\label{Xcubature} d \hat{X}_t^x(\omega_{j})&=&V_0
\bigl(\hat{X}_t^x(\omega_{j})\bigr)\,dt+\sum
_{i=1}^d V_i\bigl(
\hat{X}_t^x(\omega_{j})\bigr) \,d
\omega_{j,t}^i,
\nonumber\\[-8pt]\\[-8pt]\nonumber
\hat{X}_0^x(\omega_{j})&=&x,\qquad j=1,\ldots,l.
\end{eqnarray}
Then, for a Lipschitz continuous function $f$,
$P_tf(x)=E[f(X_t^x)]$ can be approximated by $\hat
{Q}_{(t)}^{m}f(x)=\sum_{j=1}^l \lambda_j f(\hat{X}_t^x(\omega
_{j}))$ as follows:
%
\begin{eqnarray}
\bigl\llVert P_t f - \hat{Q}_{(t)}^{m}f \bigr
\rrVert _{\infty}&=&O\bigl(t^{(m+1)/2}\bigr),\qquad t\in(0,1].
\end{eqnarray}
Then it can be shown that
%
\begin{eqnarray}
\bigl\llVert P_T f - \bigl(\hat{Q}_{(T/n)}^{m}
\bigr)^n f \bigr\rrVert _{\infty
}&=&O\bigl(n^{-(m-1)/2}
\bigr).
\end{eqnarray}
See \citeauthor{Kusuokaa} (\citeyear{Kusuokaa,Kusuokad}) and \citet{Lyons-Victoir} for the proofs.
Here, we note that the Kusuoka--Lyons--Victoir's approximation is
generally discussed in the case of nonuniform time grids.

\begin{algorithm}[t]
\caption{Weak approximation with asympotic expansion and Malliavin weights}\label{alg1}
\begin{algorithmic}
\STATE Define the Malliavin weight ${\mathcal M}^m(t,x,y)$.
\FOR{$i=1$ \TO$n$}
\STATE Simulate Gaussian random variable $\bar{X}_{T/n}^{x,(j)}$,
$j=1,\ldots,M$.
\IF{$i=1$}
\STATE$q_{n-i}(x)=\frac{1}{M}\sum_{j=1}^M f (\bar
{X}_{T/n}^{x,(j)} ){\mathcal M}^m  (T/n,x,\bar
{X}_{T/n}^{x,(j)} )$
\ELSE
\STATE$q_{n-i}(x)=\frac{1}{M}\sum_{j=1}^M q_{n-i+1} (\bar
{X}_{T/n}^{x,(j)} ){\mathcal M}^m  (T/n,x,\bar
{X}_{T/n}^{x,(j)} )$
\ENDIF
\ENDFOR
\STATE$P_Tf(x)\simeq{q}_{0}(x)$
\end{algorithmic}
\end{algorithm}
\begin{algorithm}[b]
\caption{Weak approximation: KLV cubature
on Wiener space}\label{alg2}
\begin{algorithmic}
\STATE Define the cubature paths $\omega=(\omega_1,\ldots,\omega
_l)$ and weights $\lambda=(\lambda_1,\ldots,\lambda_l)$.
\FOR{$i=1$ \TO$n$}
\STATE Solve ODE for $\hat{X}_{T/n}^{x}(\omega_j)$, $j=1,\ldots,l$.
\IF{$i=1$}
\STATE$\hat{q}_{n-i}(x)= \sum_{j=1}^l \lambda_j  f  (\hat
{X}_{T/n}^x(\omega_j)  )$
\ELSE
\STATE$\hat{q}_{n-i}(x)= \sum_{j=1}^l \lambda_j  \hat
{q}_{n-i+1} (\hat{X}_{T/n}^x(\omega_j)  )$
\ENDIF
\ENDFOR
\STATE$P_Tf(x)\simeq\hat{q}_{0}(x)$
\end{algorithmic}
\end{algorithm}

In order to obtain a local approximation, we use Malliavin's
integration by parts formula on Wiener space for a Gaussian random
variable $\bar{X}_t^{x,\varepsilon}=X_t^{0,x}+\varepsilon\frac
{\partial}{\partial\varepsilon}X_t^{x,\varepsilon}\mid _{\varepsilon=0}$.
Then the local approximation ${Q}_{(t)}^{m}f(x)$ is given by
multiplying the Malliavin weight function ${\mathcal M}^m(t,x,\bar
{X}_t^{x,\varepsilon})$ and $f(\bar{X}_t^{x,\varepsilon})$.
As mentioned above, we can mostly evaluate Malliavin weights ${\mathcal
M}^m(t,x,\bar{X}_t^{x,\varepsilon})$
as closed forms.

On the other hand, the Kusuoka--Lyons--Victoir's approximation based on
the cubature formula on Wiener space requires to solve the ODEs (\ref
{Xcubature}) with
cubature paths and weights, and then
the local approximation $\hat{Q}_{(t)}^{m}f(x)$ is given by
the weighted sum of $f(\hat{X}_t^x(\omega_{j}))$
with the cubature weights $\lambda_j$, $j=1,\ldots,l$.

Finally, we summarize the algorithms of our weak approximation method
and the KLV cubature scheme
on Wiener space
as Algorithms \ref{alg1} and~\ref{alg2}, respectively.

\section{Numerical example}\label{sec5}
This section demonstrates the effectiveness of our method
through the numerical examples for option pricing under local and
stochastic volatility models.
\subsection{Local volatility model}\label{sec5.1}
The first example takes the following local volatility model:
%
\begin{eqnarray}
dS_t^{x,\varepsilon}&=&\varepsilon\sigma\bigl(S_t^{x,\varepsilon}
\bigr) \,dB_t,
\nonumber\\[-8pt]\\[-8pt]\nonumber
S_0^{x,\varepsilon}&=&S_0 = x.
\end{eqnarray}
Then let $(\bar{S}_t^{x,\varepsilon})_{t\geq0}$ be the solution to
the following SDE:
%
\begin{eqnarray}
d\bar{S}_t^{x,\varepsilon}&=&\varepsilon\sigma(x) \,dB_t,
\nonumber\\[-8pt]\\[-8pt]\nonumber
\bar{S}_0^{x,\varepsilon}&=&x.
\end{eqnarray}
In this numerical example, for the payoff function $f(x)=\max\{x-K,0\}
$ or $f(x)=\max\{K-x,0\}$ where $K$ is a positive constant,
we apply the first-order asymptotic expansion operator, that is, $m=1$;
%
\begin{equation}
Q^1_{(t)} f(x)=E\bigl[f\bigl(\bar{S}^{x,\varepsilon}_t
\bigr){\mathcal M}^1\bigl(t,x,\bar {S}_t^{x,\varepsilon}
\bigr) \bigr]
\end{equation}
and
the second-order asymptotic expansion operator, that is, $m=2$;
%
\begin{equation}
Q^2_{(t)} f(x)=E\bigl[f\bigl(\bar{S}^{x,\varepsilon}_t
\bigr){\mathcal M}^2\bigl(t,x,\bar {S}_t^{x,\varepsilon}
\bigr) \bigr].
\end{equation}
The Malliavin weights ${{\mathcal M}}^{1}(t, x, y)$ and ${{\mathcal
M}}^{2}(t, x, y)$ are given by
\begin{eqnarray*}
&&{{\mathcal M}}^{1}(t, x, y) =1+\varepsilon E \biggl[H_{(1)}
\biggl( \frac{\partial}{\partial\varepsilon} S^{x, \varepsilon
}_{t}\bigg| _{\varepsilon=0},
\frac{1}{2}\frac{\partial^2}{\partial
\varepsilon^2} S^{x, \varepsilon}_{t}\bigg|
_{\varepsilon=0} \biggr) \Big| \bar{S}_{t}^{x,\varepsilon}=y \biggr],
\end{eqnarray*}
and
\begin{eqnarray*}
{{\mathcal M}}^{2}(t, x, y) &=&{{\mathcal M}}^{1}(t, x,
y)+\varepsilon ^2 E \biggl[H_{(1)} \biggl(
\frac{\partial}{\partial\varepsilon} S^{x, \varepsilon}_{t}\bigg| _{\varepsilon=0},
\frac{1}{6}\frac{\partial^3}{\partial\varepsilon^3} S^{x, \varepsilon}_{t} \bigg|
_{\varepsilon=0} \biggr) \Big| \bar{S}^{x}_{t}=y \biggr]
\\
&&{}+ \frac{1}{2}\varepsilon^2 E \biggl[H_{(1,1)}
\biggl(\frac{\partial
}{\partial\varepsilon} S^{x, \varepsilon}_{t}\bigg| _{\varepsilon=0},
\biggl( \frac{1}{2} \frac{\partial^2}{\partial\varepsilon^2} S^{x, \varepsilon}_{t}\bigg|
_{\varepsilon=0} \biggr)^2 \biggr) \Big| \bar{S}_{t}^{x,\varepsilon}=y
\biggr].
\end{eqnarray*}
Moreover, we remark that those Malliavin weights are obtained as closed forms.

Also, we specify the local volatility function
as a log-normal scaled volatility $\varepsilon\sigma(S) =\varepsilon
S_0^{1-\beta} S^\beta$ with $\beta=0.5$.
The parameters are set to be
$S_0=100$ and $\varepsilon= 0.4$.
The benchmark values are computed by Monte Carlo simulations
($\mathrm{Benchmark~MC}$) with
$10^7$ trials and 1000 time steps for the 1 year maturity case or 2000
time steps for the
10 year maturity case.

\begin{figure}[t]

\includegraphics{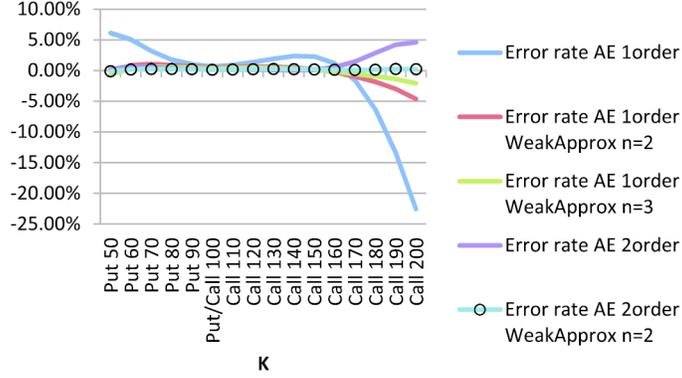}

\caption{$T=1$: Local volatility model, Error rates of the first- and
second-order asymptotic expansions and their weak approximations.}\label{fig1}
\end{figure}

\begin{figure}[b]

\includegraphics{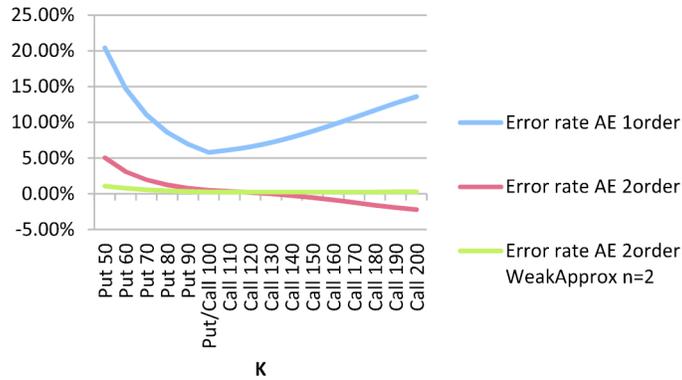}

\caption{$T=10$: Local volatility model, Error rates of the first- and
second-order asymptotic expansions and the weak approximations.}\label{fig2}
\end{figure}

Figures~\ref{fig1} and~\ref{fig2} show the results. The vertical axis in the
figures is the Error rate defined by
\begin{eqnarray*}
&&\mbox{Error Rate}=(\mathrm{WeakApprox}-\mathrm{Benchmark~MC})/\mathrm{Benchmark~MC~(\%)}.
\end{eqnarray*}
Here, $\mathrm{WeakApprox}$ is our weak approximation based on the asymptotic
expansion with Malliavin weights given in previous sections.
We observe that the increase in the number of the time steps improves
the approximation. (See Error rate AE 1order and Error rate AE 1order
WeakApprox $n=2,3$ in Figure~\ref{fig1}.)
We also note that
our scheme with the second-order expansion and two time steps
(Error rate AE 2order WeakApprox $n=2$)
improves the base (analytical only) second-order expansion (Error rate
AE 2order), and
is able to provide an accurate approximation across all the strikes
even for the long maturity case such as the 10-year maturity case in
Figure~\ref{fig2}.

\subsection{Stochastic volatility model}\label{sec5.2}
The second example considers the following stochastic volatility model,
which is also known as the log-normal SABR model:
%
\begin{eqnarray}
dS_t^{(z,\sigma)} &=&\sigma_t^{\sigma}
S_t^{(z,\sigma)} \,dB_t^1,\qquad
S_0^{(z,\sigma)} = z,
\\
d\sigma_t^{\sigma}&=&\nu\sigma_t^{\sigma}
\bigl(\rho \,dB_t^1 + \sqrt {1-\rho^2}
\,dB_t^2\bigr),\qquad \sigma_0^{\sigma}=
\sigma.
\end{eqnarray}
Next, let us introduce the following perturbed logarithmic SABR model:
%
\begin{eqnarray}
dX_{1,t}^{(x,\sigma),\varepsilon} &=&\varepsilon \biggl[-\eta
\frac{(\sigma_t^{\sigma,\varepsilon})^2}{2} \,dt +\eta\sigma _t^{\sigma,\varepsilon}
\,dB_t^1 \biggr],\qquad X_{1,0}^{(x,\sigma
),\varepsilon} = x,
\\
d\sigma_t^{\sigma,\varepsilon}&=&\varepsilon \bigl[ \sigma
_t^{\sigma,\varepsilon} \bigl(\rho \,dB_t^1 +
\sqrt{1-\rho^2} \,dB_t^2\bigr) \bigr],\qquad
\sigma_0^{\sigma,\varepsilon}=\sigma,
\end{eqnarray}
with $\varepsilon=\nu$ and $\eta=1/\nu$.
For some fixed $T>0$ and $K>0$, the target expectation is given by
\begin{eqnarray}
&& E \bigl[f\bigl(X_{1,T}^{(x,\sigma),\varepsilon}, \sigma_T^{\sigma,\varepsilon}
\bigr) \bigr]\nonumber
\\
&&\qquad \equiv E \bigl[ \hat{f}\bigl(X_{1,T}^{(x,\sigma
),\varepsilon}\bigr)
\bigr]
\nonumber
\\
&&\qquad := E \bigl[\max \bigl\{e^{X_{1,T}^{(x,\sigma),\varepsilon}} - K,0 \bigr\} \bigr]\quad\mbox{or}\quad E
\bigl[\max \bigl\{K - e^{X_{1,T}^{(x,\sigma),\varepsilon}},0 \bigr\} \bigr].
\nonumber
\end{eqnarray}

Next, let $(\bar{X}_{1,t}^{(x,\sigma),\varepsilon},\bar{\sigma
}_t^{\sigma,\varepsilon})_{t\geq0}$ be the solution to the following SDE:
%
\begin{eqnarray}
d\bar{X}_{1,t}^{(x,\sigma),\varepsilon} &=&\varepsilon \biggl[ -\eta
\frac{\sigma^2}{2} \,dt + \eta\sigma \,dB_t^1 \biggr], \qquad \bar
{X}_{1,0}^{(x,\sigma),\varepsilon} = x,
\\
d \bar{\sigma}_t^{\sigma,\varepsilon} &=&\varepsilon \bigl[ \sigma\bigl(
\rho \,dB_t^1 + \sqrt{1-\rho^2}
\,dB_t^2\bigr) \bigr],\qquad \bar {\sigma}_0^{\sigma,\varepsilon}=
\sigma.
\end{eqnarray}
The parameters are set to be $z=100$, $\sigma= 0.3$, $\varepsilon=\nu
=0.1$, $\eta=1/\nu$ and $\rho= -0.5$.
The benchmark values are calculated by Monte Carlo simulations with
$10^7$ trials and 1000 time steps for the 1-year maturity case or 2000
times steps for the 2-year maturity case.

In this example, we use the first-order two-dimensional asymptotic
expansion operator with two time steps, that is, $m=1$ and $n=2$.
Then the calculation procedure corresponding to the one in the previous
section is the following:
first, set $t_0=0$, $t_1=T/2$, $t_2=T$ and $s=t_k-t_{k-1}=T/2$ ($k=1,2$).
\begin{itemize}
\item
$\mbox{For }
 (\bar{X}_{1,t_1}^{(x_1,\sigma_1),\varepsilon},\bar{\sigma
}_{1,t_1}^{\sigma_1,\varepsilon}) = (x_1,\sigma_1) \mbox{ at } t=t_1$,
%
\begin{equation}\label{svq1}
q_1(x_1,\sigma_1)=E \bigl[\hat{f} \bigl(
\bar{X}^{(x_1,\sigma
_1),\varepsilon}_{1,s} \bigr) {\mathcal M}^1
\bigl(s,(x_1,\sigma_1), \bigl(\bar{X}^{(x_1,\sigma
_1),\varepsilon}_{1,s},
\bar{\sigma}^{\sigma_1,\varepsilon
}_{s} \bigr) \bigr) \bigr].
\end{equation}
\item
At $t=t_0=0$,
%
\begin{eqnarray}
\label{q0M1}
q_0(x,\sigma)&=&E \bigl[q_1 \bigl(
\bar{X}^{(x,\sigma),\varepsilon
}_{1,s}, \bar{\sigma}^{\sigma,\varepsilon}_{s}
\bigr) {\mathcal M}^1 \bigl(s, (x,\sigma), \bigl(
\bar{X}^{(x,\sigma
),\varepsilon}_{1,s}, \bar{\sigma}^{\sigma,\varepsilon}_{s}
\bigr) \bigr) \bigr].
\label{svq0}
\end{eqnarray}
\end{itemize}
Here, ${{\mathcal M}}^{1}(t, (x,\sigma), (x',\sigma'))$ is the
two-dimensional Malliavin weight given by
\begin{eqnarray*}
&&{{\mathcal M}}^{1}\bigl(t, (x,\sigma), \bigl(x',
\sigma'\bigr)\bigr)
\\
&&\qquad =1+\varepsilon E \biggl[H_{(1)} \biggl( \biggl(\frac{\partial}{
\partial\varepsilon}
X^{(x,\sigma), \varepsilon
}_{1,t}\bigg| _{\varepsilon=0}, \frac{\partial}{\partial\varepsilon}
\sigma^{\sigma, \varepsilon}_{t}\bigg| _{\varepsilon=0} \biggr),
\frac{1}{2}\frac{\partial^2}{\partial\varepsilon^2} X^{(x,\sigma),
\varepsilon}_{1,t}\bigg|
_{\varepsilon=0} \biggr) \Big|
\\
&&\hspace*{216pt}  \bigl(\bar {X}^{(x,\sigma),\varepsilon}_{1,t},
\bar{\sigma}^{\sigma,\varepsilon}_{t} \bigr)=\bigl(x',
\sigma'\bigr) \biggr]
\nonumber
\\
&&\quad\qquad{}+\varepsilon E \biggl[H_{(1)} \biggl( \biggl(\frac{\partial}{\partial
\varepsilon}
X^{(x,\sigma), \varepsilon}_{1,t}\bigg| _{\varepsilon=0}, \frac{\partial}{\partial\varepsilon}
\sigma^{\sigma, \varepsilon
}_{t}\bigg| _{\varepsilon=0} \biggr),
\frac{1}{2}\frac{\partial^2}{\partial
\varepsilon^2} \sigma^{\sigma, \varepsilon}_{t}\bigg|
_{\varepsilon=0} \biggr) \Big|
\\
&&\hspace*{194pt}\bigl(\bar{X}^{(x,\sigma),\varepsilon}_{1,t},
\bar{\sigma}^{\sigma,\varepsilon}_{t} \bigr)=\bigl(x',
\sigma'\bigr) \biggr].
\end{eqnarray*}
Moreover, we remark that those Malliavin weights are obtained as closed
forms as in the local volatility case.

Actually, at $t_1$ we need not implement (\ref{svq1}),
but just compute the first-order analytical asymptotic expansion for
pricing\vspace*{1pt} options
with the time-to-maturity $T/2$ and the initial value $(\bar
{X}_{1,t_1}^{(x_1,\sigma_1),\varepsilon},\bar{\sigma}_{t_1}^{\sigma
_1,0}) = (x_1,\sigma_1)$.
That is,
%
\begin{equation}
\label{svM1} \hat{q}_1(x_1,\sigma_1)=E
\bigl[\hat{f} \bigl(\bar {X}^{(x_1,\sigma_1),\varepsilon}_{1,s} \bigr) \hat{{\mathcal
M}}^1 \bigl(s,(x_1,\sigma_1),
\bar{X}^{(x_1,\sigma
_1),\varepsilon}_{1,s} \bigr) \bigr], \label{svq10}
\end{equation}
where $\hat{{\mathcal M}}^1(s, (x_1,\sigma_1), y) = 1+\varepsilon
\hat{{\mathcal M}}_{(1)}(s, (x_1,\sigma_1), y)$, and
$\hat{{\mathcal M}}_{(1)}(s, (x_1,\sigma_1), y)$ stands for the
first-order one-dimensional Malliavin weight:
%
\begin{eqnarray}
\qquad && \hat{{\mathcal M}}_{(1)}\bigl(s, (x_1,
\sigma_1), y\bigr)
\nonumber\\[-8pt]\\[-8pt]\nonumber
&&\qquad = E \biggl[H_{(1)} \biggl(
\frac{\partial}{\partial\varepsilon} X^{(x_1,\sigma_1), \varepsilon}_{1,s}\bigg| _{\varepsilon=0},
\frac{1}{2}\frac{\partial^2}{\partial\varepsilon^2} X^{(x_1,\sigma
_1), \varepsilon}_{1,s}\bigg|
_{\varepsilon=0} \biggr) \Big| \bar{X}^{(x_1,\sigma_1,\varepsilon)}_{1,s}=y \biggr].
\end{eqnarray}
On the other hand, we apply a conditional expectation formula for
multidimensional asymptotic expansions
in \citet{Takahashi} in order to evaluate the Malliavin weight
${\mathcal M}^1$ in (\ref{q0M1}).

\begin{figure}[b]

\includegraphics{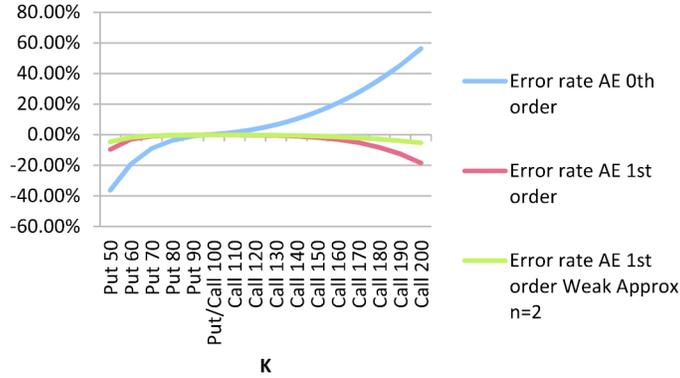}

\caption{$T=1$: Stochastic volatility model, Error rate of the
first-order two-dimensional asymptotic expansion and the weak approximation.}\label{fig3}
\end{figure}

\begin{figure}[t]

\includegraphics{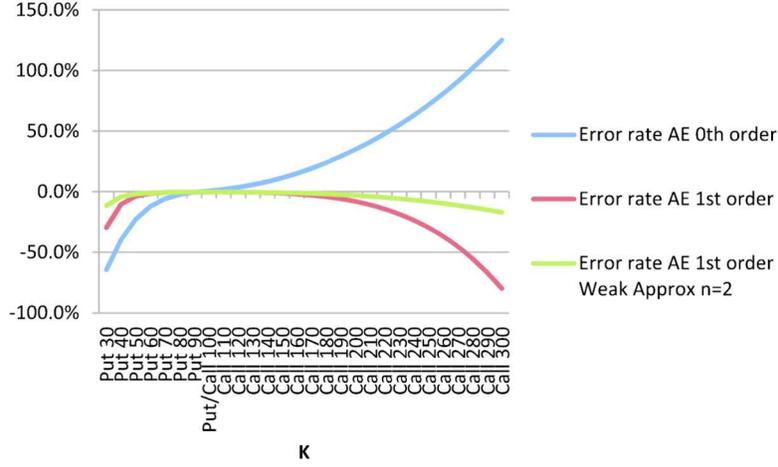}

\caption{$T=2$: Stochastic volatility model, Error rate of the
first-order two-dimensional asymptotic expansion and the weak approximation.}\label{fig4}
\end{figure}

Figures~\ref{fig3} and~\ref{fig4} show the results (the vertical axis in the
figures is Error rate).
Again, our scheme with (\ref{svM1})
and (\ref{svq0}) (Error rate AE 1st order WeakApprox $n=2$)
improves the base first-order expansion (Error rate AE 1st order)
especially for the deep OTM calls and puts.

\subsection{Error analysis}\label{sec5.3}
In this section, we investigate
the validity of our approximation by comparing
the theoretical and the numerical errors.

In particular, we use the same example of our weak approximation in the
local volatility model with the maturity $T=1$ in Section~\ref{sec5.1}. Here, we
remark that this example can be regarded as the Lipschitz continuous
case in Theorem~\ref{33333333333333}.

Based on the results of Theorem~\ref{33333333333333}, for a fixed expansion order $m$ and
a time grid parameter $\gamma=1$, the error of the weak approximation
by the $m$th order asymptotic expansion with discretization $n$
approximately satisfies the following relation:
%
\begin{eqnarray}\label{TheoreticalErrorFormula1}
&&\mbox{Error Weak Approx with Asymp Expansion } (m,n+1)\nonumber
\\
&&\qquad \simeq \bigl\{\mbox{Error Weak Approx with Asymp Expansion }(m,n) \bigr\}
\\
&&\quad\qquad{}\times \bigl(n/(n+1)\bigr)^{m/2}.\nonumber
\end{eqnarray}
Here, ``\textit{Error Weak Approx with Asymp Expansion} ($m,n$)''
stands for
the deviation of ``Weak Approximation'' from ``Benchmark Monte Carlo,''
that is
the value of $(\mathrm{WeakApprox})-(\mathrm{Benchmark~MC})$.

\begin{figure}

\includegraphics{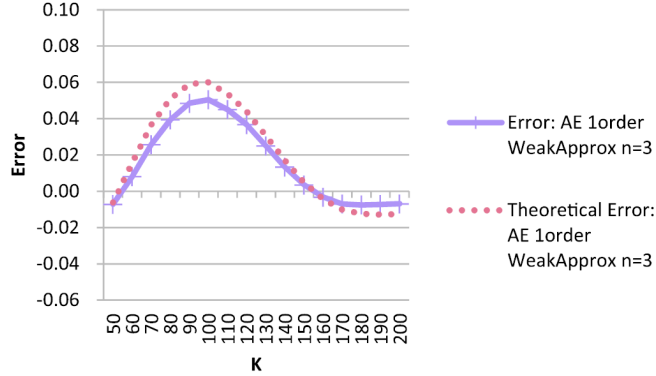}

\caption{Error with respect to $n$ of the weak approximation for fixed $m$.}\label{fig5}
\end{figure}

Figure~\ref{fig5} checks the above relation in the case that $m=1$ and $n=2$.
In the figure,
``\textit{Theoretical Error}: \textit{AE} 1\textit{order WeakApprox} $n=3$'' is calculated by
equation~(\ref{TheoreticalErrorFormula1}). It is observed that
the order of the theoretical error is rather
similar to that of the numerical error ``\textit{Error}: \textit{AE} 1\textit{order WeakApprox} $n=3$''
across all the strike prices.

\begin{figure}[b]

\includegraphics{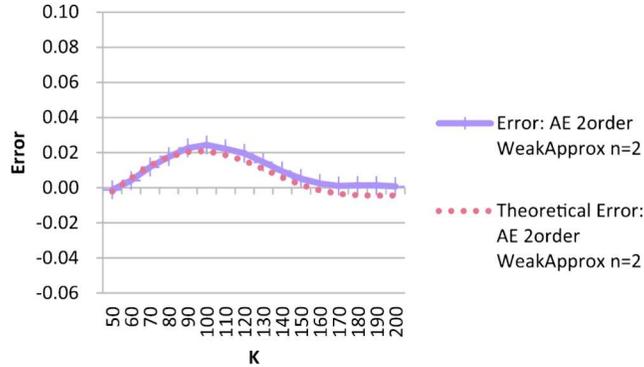}

\caption{Error with respect to $m$ of the weak approximation for fixed $n$.}\label{fig6}
\end{figure}

Next, let us check the validity of our method
from another viewpoint.
For a fixed partition number $n$ and a time grid parameter $\gamma=1$,
the error of the weak approximation based on ($m+1$)th order
asymptotic expansion with discretization $n$
labeled by
``\textit{Error Weak Approx with Asymp Expansion} ($m+1,n$),'' approximately
satisfies the relation:
\begin{eqnarray}\label{TheoreticalErrorFormula2}
\qquad &&\mbox{Error Weak Approx with Asymp Expansion }(m+1,n)\nonumber
\nonumber\\[-8pt]\\[-8pt]\nonumber
&&\qquad \simeq \bigl\{\mbox{Error Weak Approx with Asymp Expansion }(m,n)\bigr\} \times
\varepsilon/ \sqrt{n}.
\end{eqnarray}
Figure~\ref{fig6} examines
the above relation in the case $m=1$ and $n=2$ with $\varepsilon=0.4$.

``\textit{Theoretical Error}: \textit{AE} 2\textit{order WeakApprox} $n=2$'' in Figure~\ref{fig6} is calculated
by the equation (\ref{TheoreticalErrorFormula2}). We observe
that the order of the theoretical error is very close to that of the
numerical error
``\textit{Error}: \textit{AE} 2\textit{order WeakApprox} $n=2$'' for all the strike prices.

Finally,
we test the numerical errors by changing the parameter $\gamma$.
Again, we fix the parameter $m=1$.
Based on the result of the Lipschitz continuous case in our main
theorem (Theorem~\ref{33333333333333}), the errors depend on the range of
$\gamma$, {that is}, $\gamma<1/3=m/(m+2)$, $\gamma=1/3=m/(m+2)$,
$\gamma>1/3=m/(m+2)$.

\begin{figure}[t]

\includegraphics{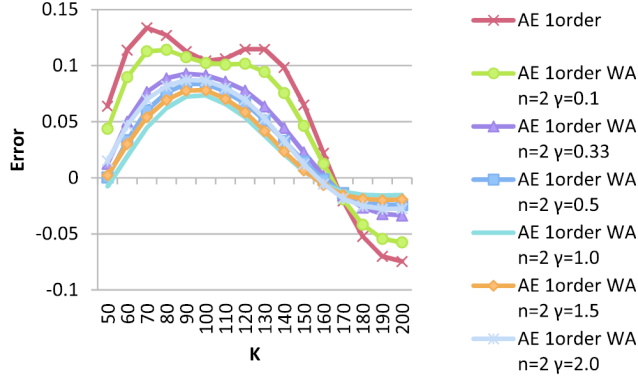}

\caption{Error of the weak approximation $n=2$ with respect to time
grid parameter $\gamma$.}\label{fig7}
\end{figure}

\begin{figure}[b]

\includegraphics{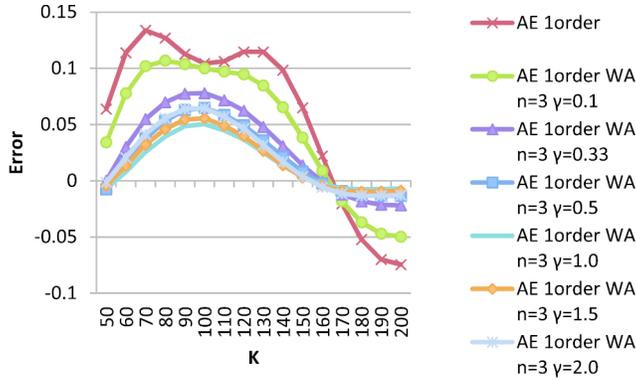}

\caption{Error of the weak approximation $n=3$ with respect to time
grid parameter $\gamma$.}\label{fig8}
\end{figure}

In order to see the differences of the errors with the different values
of $\gamma$,
Figures~\ref{fig7} and~\ref{fig8} plot
the errors for $\gamma=0.1$, $\gamma=0.33$, $\gamma=0.5$, $\gamma
=1.0$, $\gamma=1.5$ and $\gamma=2.0$ with $n=2$ and $n=3$, respectively.

We are able to find that the errors are determined by the levels of
$\gamma$ and
the behavior of the errors is consistent with the theoretical results
in Theorem~\ref{33333333333333}.

In addition, we examine which time grid parameter $\gamma$ is optimal.
In order to show this, we execute a simple test for the case $m=1$ and
$n=2$. Particularly, we solve the following minimization problem:
%
\begin{eqnarray}
\hat{\gamma} &=&\operatorname{argmin} \biggl\{ \sum_{K \in\{50,60,\ldots,190,200 \}}
\nonumber\\[-8pt]\\[-8pt]\nonumber
&&\hspace*{38pt} \bigl(\mbox{Error Weak Approx with Asymp Expansion }(m,n; \gamma,K)
\bigr)^2 \biggr\}.\hspace*{-20pt}
\end{eqnarray}

We obtained a parameter $\hat{\gamma}=1.015657$.
That is, the value close to
$\gamma=1$ (the uniform time grid case) is optimal in our weak
approximation with asymptotic expansion.

Therefore, we can conclude that the results of the numerical
experiments of our weak approximation are consistent with the
theoretical part of this paper, and confirmed the validity of our method.

\section{Concluding remarks}\label{sec6}
In this paper, we have shown a new approximation method for the
expectations of the functions of the solutions to SDEs
by applying an asymptotic expansion with Malliavin calculus. In
particular, based on \citeauthor{Kusuokaa}
(\citeyear{Kusuokaa,Kusuokab,Kusuokac,Kusuokad}),
we have obtained error estimates for our new weak approximation.

Moreover, we have confirmed the validity of our method
through the numerical examples for option pricing under local and
stochastic volatility models.
The scheme is simple and we can attain enough accuracy even when the
expansion order $m$ is low such as $m=1,2$ with a few time steps
$n=2,3$ as demonstrated in the previous section.

In order to obtain more accurate numerical approximation,
it is natural to use many partitions $n$ in the time scale.
However, the computational cost becomes exponentially larger as the
number of
partitions becomes larger.

To overcome this problem, some efficient tree based (discretization) techniques
can be applied.
Another possible solution is to use the higher order expansion
developed in \citet{Takahashi-Takehara-Todab} or \citet{Violante}. We are
convinced that the higher order expansion will improve the accuracy
since the higher order $m$th expansion improves the error orders to
$O(\varepsilon^{m+1}/n^{m})$ for a Lipschitz continuous $f$ and
$O(\varepsilon^{m+1}/n^{m-1})$
for a bounded Borel $f$.

Further, applying our method to
the higher-dimensional problems is one of the important issues. When
the dimension $N$ of the state variables becomes higher, the
computational cost becomes larger.
However,
the multidimensional higher order expansion such as in \citet{Takahashi} or
\citet{Takahashi-Takehara-Todab}
is a tractable approach to the extension.
These topics will be the main themes in our next research.

\begin{appendix}\label{appe}
\section{Proof of Theorem~\texorpdfstring{\protect\ref{11111111111111}}{1}}\label{appeA}
First, for the preparation for the proof of the theorem, we
characterize the differentiations of the solution to the general
perturbed SDEs $X_t^{x,\varepsilon}$ with respect to $\varepsilon$
as elements in the space ${\mathcal K}_r$. The following lemma plays an
important rule for estimating the order of the local approximation for
$E[f(X_t^{x,\varepsilon})]$ in Theorem~\ref{11111111111111}.

\begin{lemma}\label{differentiationKr}\label{le2}
\begin{eqnarray*}
&&\frac{1}{j!} \frac{\partial^j}{\partial\varepsilon^j} X_t^{x,\varepsilon} \in {
\mathcal K}_j,\qquad j \geq1.
\end{eqnarray*}
\end{lemma}

\begin{pf}
We prove the assertion by induction.
First, the differentiation of $X_{t}^{x,\varepsilon}$ with respect to
$\varepsilon$ is given by
%
\begin{eqnarray}\label{1stSDE}
\frac{\partial}{\partial\varepsilon}X_{t}^{x,\varepsilon,l} &=&\int_{0}^{t}
\frac{\partial}{\partial\varepsilon} V_{0}^l\bigl(\varepsilon,
X_{s}^{x,\varepsilon}\bigr)\,ds+\sum_{j=1}^{d}
\int_{0}^{t}V_{j}^l
\bigl(X_{s}^{x,\varepsilon}\bigr)\,dB_{s}^{j}
\nonumber
\\
&&{}+ \sum_{k=1}^N \int
_{0}^{t} \partial_k
V_{0}^l\bigl(\varepsilon,X_{s}^{x,\varepsilon}
\bigr)\frac{\partial}{\partial\varepsilon
}X_{s}^{x,\varepsilon,k}\,ds
\\
&&{}+ \varepsilon \sum_{k=1}^N \sum
_{j=1}^{d}\int_{0}^{t}
\partial_k V_{j}^l\bigl(X_{s}^{x,\varepsilon}
\bigr)\frac{\partial}{\partial\varepsilon
}X_{s}^{x,\varepsilon,k}\,dB_{s}^{j},\qquad
l=1,\ldots,N.\nonumber
\end{eqnarray}
The above SDE is linear and
the order of the Kusuoka--Stroock function $\frac{\partial}{\partial
\varepsilon}X_{t}^{x,\varepsilon}$ is determined by the following term:
%
\begin{equation}
\sum_{j=1}^d \int_0^t
J_t^{x,\varepsilon} \bigl(J_u^{x,\varepsilon}
\bigr)^{-1} V_j\bigl({X}_u^{x,\varepsilon}
\bigr) \,dB_u^j \in{\mathcal K}_1,
\end{equation}
where $J_t^{x,\varepsilon}=\nabla_x X_t^{x,\varepsilon}$.
Since this term gives the minimum order in the terms that consist of
(\ref{1stSDE}). Here, we use the properties $J_s^{x,\varepsilon},
(J_s^{x,\varepsilon})^{-1} \in{\mathcal K}_0$, $s \in(0,1]$ and the
boundness of $V_j$, $j=1,\ldots,d$. We have $\frac{\partial
}{\partial\varepsilon}X_{t}^{x,\varepsilon} \in{\mathcal K}_1$ by
using the properties 2 and 3 in Lemma \ref{Properties}.

For $i\geq2$, $\frac{1}{i!}\frac{\partial^i}{\partial\varepsilon
^i}X_{t}^{x,\varepsilon}= ( \frac{1}{i!}\frac{\partial
^i}{\partial\varepsilon^i}X_{t}^{x,\varepsilon,1},\ldots, \frac
{1}{i!}\frac{\partial^i}{\partial\varepsilon^i}X_{t}^{x,\varepsilon,N}  )$ is recursively determined by the following:
%
\begin{eqnarray}\label{ExpansionSDE}
&& \frac{1}{i!}\frac{\partial^i}{\partial\varepsilon
^i}X_{t}^{x,\varepsilon,n}\nonumber
\\
&&\qquad =
\frac{1}{i!}\int_0^t
\frac{\partial^i}{\partial\varepsilon
^i}V_0^n\bigl(\varepsilon,
X_u^{x,\varepsilon}\bigr)\,du\nonumber
\\
&&\quad\qquad{}+\sum_{m=1}^i \sum
_{{i}^{(k)},{\alpha}^{(k)}}^{(m)} \frac{1}{(i-m)!} \int
_0^t \Biggl( \prod
_{l=1}^{k} \frac{1}{i_l !} \frac{\partial
^{i_l}}{\partial\varepsilon^{i_l}}X_{u}^{x,\varepsilon, \alpha_l}
\Biggr) \partial_{{\alpha}^{(k)}} \frac{\partial^{i-m}}{\partial
\varepsilon^{i-m}} V_0^n
\bigl(\varepsilon,{X}_u^{x,\varepsilon}\bigr)\,du\hspace*{-25pt}
\nonumber\\[-8pt]\\[-8pt]\nonumber
&&\quad\qquad{}+\sum_{{i}^{(k)},{\alpha}^{(k)}}^{(i-1)} \int
_0^t \Biggl( \prod
_{l=1}^{k} \frac{1}{i_l !} \frac{\partial
^{i_l}}{\partial\varepsilon^{i_l}}X_{u}^{x,\varepsilon, \alpha_l}
\Biggr)\sum_{j=1}^d \partial_{{\alpha}^{(k)}}
V_j^n\bigl({X}_u^{x,\varepsilon}
\bigr)\,dB_u^j
\nonumber
\\
&&\quad\qquad{}+\varepsilon\sum_{{i}^{(\beta)},{\alpha}^{(k)}}^{(i)} \int
_0^t \Biggl( \prod
_{l=1}^{k} \frac{1}{i_l !} \frac{\partial
^{i_l}}{\partial\varepsilon^{i_l}}X_{u}^{x,\varepsilon, \alpha_l}
\Biggr) \sum_{j=1}^d
\partial_{{\alpha}^{(k)}} V_j^n\bigl({X}_u^{x,\varepsilon}
\bigr)\,dB_u^j,\nonumber
\\[-4pt]
\eqntext{n=1,\ldots,N,}
\end{eqnarray}
where
%
\begin{equation}
\sum_{{i}^{(k)},{\alpha}^{(k)}}^{(i)}:= \sum
_{k=1}^i \sum_{i_1+\cdots+i_k=i, i_l \geq1}
\sum_{{\alpha
^{(k)}} \in\{1,\ldots,N\}^k} \frac{1}{k !}.
\end{equation}
The above SDE is linear and the order of the Kusuoka--Stroock function
$\frac{1}{i!}\frac{\partial^i}{\partial\varepsilon
^i}X_{t}^{x,\varepsilon}$ is determined inductively by the term
%
\begin{eqnarray}
\qquad &&\sum_{{i}^{(k)},{\alpha}^{(k)}}^{(i-1)} \int
_0^t J_t^{x,\varepsilon}
\bigl(J_u^{x,\varepsilon} \bigr)^{-1} \Biggl( \prod
_{l=1}^{k} \frac{1}{i_l !} \frac{\partial^{i_l}}{\partial\varepsilon
^{i_l}}X_{u}^{x,\varepsilon, \alpha_l}
\Biggr)\sum_{j=1}^d \partial
_{{\alpha}^{(k)}} V_j\bigl(X_u^{x,\varepsilon}
\bigr)\,dB_s^j \in{\mathcal K}_i.
\end{eqnarray}
Since this term gives the minimum order in the terms that consist of
(\ref{ExpansionSDE}).
Then $\frac{1}{i!}\frac{\partial^i}{\partial\varepsilon
^i}X_{t}^{x,\varepsilon} \in{\mathcal K}_i$ by using the properties 2
and 3 in Lemma \ref{Properties}.
\end{pf}

Hereafter, we give the expansion for $E[f(X_t^{x,\varepsilon})]$
around $E[f(\bar{X}_t^{x,\varepsilon})]$. We remark that
$X_t^{x,\varepsilon}$ is not uniformly nondegenerate Wiener functional
in Watanabe sense because $X_t^{x,0}$ is completely degenerate as
Wiener functional, {that is}, $X_t^{x,0}$ is the solution to ODE.
Then, in order to give the expansion, we define a Wiener functional
$Y_t^{\varepsilon}$ given by $Y_t^{\varepsilon}=\varphi
(X_t^{x,\varepsilon})=\frac{X_t^{x,\varepsilon
}-X_t^{x,0}}{\varepsilon}$, {that is}, $(Y_t^{\varepsilon,1},\ldots,Y_t^{\varepsilon,N})=(\varphi_1(X_t^{x,\varepsilon,1}),\ldots,\varphi_N(X_t^{x,\varepsilon,N}))$, $\varphi_i(\xi
)=\frac{\xi-X_t^{x,0,i}}{\varepsilon}$, $i=1,\ldots,N$. The
expansion of $Y_t^{\varepsilon}$ is given in the space ${\mathbf
D}^{\infty}$, that is, for all $m \in{\mathbf N}$,
%
\begin{eqnarray}
\qquad && \limsup_{\varepsilon\downarrow0} \frac{1}{\varepsilon^{m+1}} \Biggl\llVert
Y_t^{\varepsilon}- \Biggl\{ \frac{\partial}{\partial
\varepsilon} X_{t}^{x,\varepsilon}
\bigg| _{\varepsilon=0}+\sum_{i=1}^m
\varepsilon^i \frac{1}{(i+1)!}\frac{\partial^{i+1}}{\partial
\varepsilon^{i+1}}X_t^{x,\varepsilon}
\bigg| _{\varepsilon=0} \Biggr\} \Biggr\rrVert _{{\mathbf D}^{k,p}}
\nonumber\\[-8pt]\\[-8pt]\nonumber
&&\qquad < \infty,\qquad \forall k \in{\mathbf N}, \forall p < \infty.
\end{eqnarray}
We note that $ Y_t^{0}=\frac{\partial}{\partial\varepsilon}
X_{t}^{x,\varepsilon}\mid _{\varepsilon=0}$ and $Y_0^{0}=0$.
Let $\sigma^{Y_t^{\varepsilon}}$ be the Malliavin covariance matrix
of $Y_t^{\varepsilon}$ and
set
%
\begin{equation}
\tau=\inf \bigl\{ s; \bigl(J_s^{x,\varepsilon}\bigr)^{-1}A
\bigl(X_s^{x,\varepsilon
}\bigr) \bigl(\bigl(J_s^{x,\varepsilon}
\bigr)^{-1}\bigr)^{\top} \leq A(x)/2 \bigr\}.
\end{equation}
Then we can see
%
\begin{eqnarray}
\det\bigl(\sigma^{Y_t^{\varepsilon}}\bigr) &\geq& \det\bigl(J_t^{x,\varepsilon}
\bigr)^2 \det\int_0^{\min\{t,\tau\}}
\bigl(J_s^{x,\varepsilon
}\bigr)^{-1}A\bigl(X_s^{x,\varepsilon}
\bigr) \bigl(\bigl(J_s^{x,\varepsilon}\bigr)^{-1}
\bigr)^{\top} \,ds
\nonumber\\[-8pt]\label{A1} \\[-8pt]\nonumber
&\geq& \bigl(1/2^N\bigr)\det\bigl(J_t^{x,\varepsilon}
\bigr)^2 \det\bigl(A(x)\bigr)\min\{t,\tau\} ^N,
\\
\qquad &&\hspace*{-52pt} \sup_{\varepsilon\in(0,1]} \bigl\llVert \det\bigl(J_t^{x,\varepsilon}
\bigr)^{-1} \bigr\rrVert _{L^p} < \infty, \label{A2}
\end{eqnarray}
and
%
\begin{equation}
P(\tau<1/n) \leq c_1 \exp\bigl( -c_2 n^{c_3}
\bigr),\qquad n \in{\mathbf N}, \label{A3}
\end{equation}
where $c_i$, $i=1,2,3$ are positive constants [see the proofs of
Theorem~3.4 of \citet{Watanabe} or Theorem 10.5 of \citet{Ikeda-Watanabe} for (\ref{A1}), (\ref{A2}) and (\ref{A3})].

Therefore, under condition \textup{[\textbf{H}]}, we can see the nondegeneracy of
the Malliavin covariance matrix of $Y_t^{\varepsilon}$
%
\begin{equation}
\sup_{\varepsilon\in(0,1]} \bigl\llVert \det\bigl(\sigma^{Y_t^{\varepsilon
}}
\bigr)^{-1} \bigr\rrVert _{L^p}<\infty,\qquad p <\infty.
\end{equation}
Then the density $\xi\mapsto p_t^{Y^{\varepsilon}}(\xi)$ of
$Y_t^{\varepsilon}$ starting from $0$ is smooth. Moreover, the
Malliavin covariance matrix $\sigma^{Y_t^{\varepsilon}}$ is nondegenerate
uniformly in $\varepsilon$:
%
\begin{eqnarray}
&& \limsup_{\varepsilon\downarrow0} \bigl\llVert \det\bigl(\sigma
^{Y_t^{\varepsilon}}\bigr)^{-1} \bigr\rrVert _{L^p}=\bigl\llVert
\det\bigl(\sigma^{Y_t^{0}}\bigr)^{-1} \bigr\rrVert
_{L^p}<\infty,\qquad p <\infty.
\end{eqnarray}
Then we are able to give the following Taylor formulas for $\xi\mapsto
p_t^{Y^{\varepsilon}}(\xi)$ and $E[f(Y_t^{\varepsilon})]$ using the
Malliavin weights:
\begin{eqnarray}
p_t^{Y^{\varepsilon}}(\xi) &=& p_t^{Y^{0}}(\xi)
+\sum_{j=1}^m \varepsilon^j E\bigl[\Phi_t^j\mid Y_t^{0}=
\xi\bigr]p_t^{Y^{0}}(\xi ) \nonumber
\\
&&{}+\varepsilon^{m+1} \int_0^1(1-u)^m (m+1)
\nonumber\\[-8pt]\\[-8pt]\nonumber
&&\hspace*{10pt}{}\times  \sum_{\alpha^{(k)},\beta^{(k)}}^{m+1} E \Biggl[
H_{\alpha^{(k)}} \Biggl(Y_t^{\varepsilon u},\prod
_{l=1}^k \frac{1}{\beta_l !} \frac{\partial^{\beta_{l}}}{\partial\eta
^{\beta_{l}}}X_t^{x,\eta,\alpha_{l}}\bigg| _{\eta=\varepsilon u} \Biggr) \Big|
Y_t^{\varepsilon u}=\xi \Biggr]\nonumber
\\
&&\hspace*{10pt}{}\times  p_t^{Y^{\varepsilon u}}(\xi)\,du,\nonumber
\\
E\bigl[f\bigl(Y_t^{\varepsilon}\bigr)\bigr]&=&\int
_{{\mathbf R}^N}f(\xi)p_t^{Y^{\varepsilon
}}(\xi)\,d\xi\nonumber
\\
&=&\int_{{\mathbf R}^N}f(\xi)p_t^{Y^{0}}(\xi)\,d\xi
+\sum_{j=1}^m \varepsilon^j
\int_{{\mathbf R}^N}f(\xi)E\bigl[\Phi _t^j
\mid Y_t^{0}=\xi\bigr]p_t^{Y^{0}}(
\xi)\,d\xi \nonumber
\\
&&{}+\varepsilon^{m+1} \int_0^1(1-u)^m (m+1)
\nonumber\\[-8pt]\label{ExpansionY} \\[-8pt]\nonumber
&&\hspace*{10pt}{}\times  \sum_{\alpha^{(k)},\beta^{(k)}}^{m+1} \int
_{{\mathbf R}^N}f(\xi) E \Biggl[ H_{\alpha^{(k)}} \Biggl(Y_t^{\varepsilon u},
\prod_{l=1}^k \frac{1}{\beta_l !}\frac
{\partial^{\beta_{l}}}{\partial\eta^{\beta_{l}}} X_t^{x,\eta,\alpha_{l}}\bigg| _{\eta=\varepsilon u} \Biggr)\Big|\nonumber
\\
&&\hspace*{195pt}  Y_t^{\varepsilon
u}=\xi \Biggr]p_t^{Y^{\varepsilon u}}(
\xi)\,d\xi \,du
\nonumber
\\
&=& E\bigl[f\bigl(Y_t^{0}\bigr)\bigr]+ \sum
_{j=1}^m \varepsilon^j E\bigl[f
\bigl(Y_t^{0}\bigr)\Phi_t^j
\bigr]+\varepsilon^{m+1} r_m(t,x,\varepsilon).
\nonumber
\end{eqnarray}
Here, $\Phi_t^j$ is the Malliavin weight given by
%
\begin{eqnarray}
\Phi_t^j&=&\sum_{\alpha^{(k)},\beta^{(k)}}^j
H_{{\alpha}^{(k)}} \Biggl(Y_t^0, \prod
_{l=1}^k \frac{1}{\beta_l !}\frac{\partial^{\beta
_{l}}}{\partial\varepsilon^{\beta_{l}}}
X_t^{x,\varepsilon,\alpha
_{l}}\bigg| _{\varepsilon=0} \Biggr),
\end{eqnarray}
with
%
\begin{equation}
\sum_{\alpha^{(k)},\beta^{(k)}}^j= \sum
_{k=1}^j \sum_{\sum_{l=1}^k \beta_l =j+k,  \beta_l \geq2}
\sum_{\alpha^{(k)} = (\alpha_1,\ldots, \alpha_k) \in\{1,\ldots,N
\}^k} \frac{1}{k!}
\end{equation}
and $r_m(t,x,\varepsilon)$ is the residual:
\begin{longlist}
\item[1.]
%
\begin{eqnarray}\label{residualY1}
&& r_m(t,x,\varepsilon)\nonumber
\\
&&\qquad =\int_0^1(1-u)^m (m+1)
\\
&&\hspace*{45pt}{}\times  \sum
_{\alpha^{(k)},\beta^{(k)}}^{m+1} E \Biggl[ \partial_{\alpha^{(k)}}
f\bigl(Y_t^{\varepsilon u}\bigr) \prod_{l=1}^k
\frac{1}{\beta_l !} \frac{\partial^{\beta_{l}}}{\partial
\eta^{\beta_{l}}} X_t^{x,\eta,\alpha_{l}}\bigg|
_{\eta=\varepsilon u} \Biggr]\,du\nonumber
\end{eqnarray}
for $f \in C_b^{\infty}({\mathbf R}^N)$,
\item[2.]
\begin{eqnarray}\label{residualY2}
&& r_m(t,x,\varepsilon)\nonumber
\\
&&\qquad =\int_0^1(1-u)^m (m+1)\nonumber\\[-8pt]\\[-8pt]\nonumber
&&\hspace*{45pt}{}\times  \sum_{\alpha^{(k)},\beta^{(k)}}^{m+1} E \Biggl[
\partial_{\alpha^{(1)}} f\bigl(Y_t^{\varepsilon u}\bigr)
H_{\alpha
^{(k-1)}} \Biggl(Y_t^{\varepsilon u},\prod
_{l=1}^k \frac{1}{\beta_l
!} \frac{\partial^{\beta_{l}}}{\partial\eta^{\beta_{l}}}
X_t^{x,\eta,\alpha_{l}}\bigg| _{\eta=\varepsilon u} \Biggr) \Biggr]\,du\nonumber\hspace*{-3pt}
\end{eqnarray}
for $f \in C_b^{1}({\mathbf R}^N)$,
\item[3.]
\begin{eqnarray}\label{residualY3}
\qquad &&r_m(t,x,\varepsilon)\nonumber
\\
&&\qquad =\int_0^1(1-u)^m (m+1)
\nonumber\\[-8pt]\\[-8pt]\nonumber
&&\hspace*{45pt}{}\times  \sum_{\alpha^{(k)},\beta^{(k)}}^{m+1} E \Biggl[ f
\bigl(Y_t^{\varepsilon u}\bigr) H_{\alpha^{(k)}}
\Biggl(Y_t^{\varepsilon u},\prod_{l=1}^k
\frac{1}{\beta_l !} \frac
{\partial^{\beta_{l}}}{\partial\eta^{\beta_{l}}} X_t^{x,\eta,\alpha_{l}}\bigg|
_{\eta=\varepsilon u} \Biggr) \Biggr]\,du\nonumber
\end{eqnarray}
for an arbitrary bounded continuous function $f$.
\end{longlist}

Then, by the transformation $X_t^{x,\varepsilon}=X_t^{x,0}+\varepsilon
Y_t^{\varepsilon}$, the density $y \mapsto p^{X^{\varepsilon
}}(t,x,y)$ of $X_t^{x,\varepsilon}$ is given by
%
\begin{eqnarray}
p^{X^{\varepsilon}}(t,x,y)&=&p^{Y^{\varepsilon}}_t \bigl( \varphi (y)
\bigr)\det\biggl\llvert \frac{\partial( \varphi_1,\ldots,\varphi
_N )}{\partial(y_1,\ldots,y_N)} \biggr\rrvert
\\
&=&p^{Y^{\varepsilon}}_t \bigl(\bigl(y-X_t^{x,0}
\bigr)/\varepsilon \bigr)\frac{1}{\varepsilon^{N}}.
\end{eqnarray}
Here, we note that
\begin{eqnarray}
&&\int_{{\mathbf R}^N}f(y)p_t^{Y^{0}} \bigl(
\bigl(y-X_t^{x,0}\bigr)/{\varepsilon} \bigr)
\frac{1}{\varepsilon^{N}}\,dy\nonumber
\\
&&\qquad =\int_{{\mathbf R}^N}f(y) \frac{1}{ (2\pi\varepsilon^2)^{N/2} \det(\Sigma(t))^{1/2}}
\nonumber\\[-8pt]\\[-8pt]\nonumber
&&\quad\qquad{}\hspace*{17pt}\times  e^{-\afrac{(y-\varepsilon\mu(t)-X_t^{x,0})^{\top} \Sigma^{-1}(t)
(y-\varepsilon\mu(t)-X_t^{x,0}) }{2 \varepsilon^2}}\,dy
\nonumber
\\
&&\qquad =\int_{{\mathbf R}^N}f(y) p^{\bar{X}^\varepsilon}(t,x,y) \,dy=E\bigl[f\bigl(
\bar {X}_t^{x,\varepsilon}\bigr)\bigr],\nonumber
\end{eqnarray}
where $\mu(t)$ and $\Sigma(t)$ are the mean and the covariance matrix
of $Y_t^{0}$
and $y \mapsto p^{\bar{X}^{\varepsilon}}(t,x,y)$ is the density of
$\bar{X}_t^{x,\varepsilon}$. Also, for $G(t,x) \in{\mathcal K}_r$,
we have
%
\begin{eqnarray}
&&\int_{{\mathbf R}^N}f(y)E \bigl[H_{(i)}
\bigl(Y^{0}_t,G(t,x)\bigr) \mid Y^{0}_t=
\bigl(y-X_t^{x,0}\bigr)/{\varepsilon} \bigr]\nonumber
\\
&&\quad {}\times p_t^{Y^{0}}
\bigl(\bigl(y-X_t^{x,0}\bigr)/{\varepsilon} \bigr)
\frac{1}{\varepsilon
^{N}}\,dy
\nonumber\\[-8pt]\\[-8pt]\nonumber
&&\qquad = \int_{{\mathbf R}^N}f(y)E \bigl[H_{(i)}
\bigl(Y^{0}_t,G(t,x)\bigr) \mid \bar {X}^{x,\varepsilon}_t=y
\bigr]p^{\bar{X}^\varepsilon
}(t,x,y)\,dy
\\
&&\qquad = E\bigl[f\bigl(\bar{X}^{x,\varepsilon}_t\bigr)H_{(i)}
\bigl(Y^{0}_t,G(t,x)\bigr)\bigr],
\nonumber
\end{eqnarray}
and
%
\begin{eqnarray}
&&\int_{{\mathbf R}^N}f(y)E \bigl[H_{(i)}
\bigl(Y^{\varepsilon u}_t,G(t,x)\bigr) \mid Y^{\varepsilon u}_t=
\bigl(y-X_t^{x,0}\bigr)/{\varepsilon} \bigr]\nonumber
\\
&&\quad{}\times p_t^{Y^{\varepsilon u}}
\bigl(\bigl(y-X_t^{x,0}\bigr)/{\varepsilon} \bigr)
\frac{1}{\varepsilon^{N}}\,dy
\\
&&\qquad = E\bigl[f\bigl(\tilde{X}^{x,\varepsilon u}_t\bigr)H_{(i)}
\bigl(Y^{\varepsilon
u}_t,G(t,x)\bigr)\bigr],
\nonumber
\end{eqnarray}
with $\tilde{X}^{x,\varepsilon u}_t=X_t^{x,0}+\varepsilon
Y^{\varepsilon u}_t$, $u \in[0,1]$.

Therefore, (\ref{ExpansionY}) with (\ref{residualY1}), (\ref
{residualY2}) and (\ref{residualY3}) can be transformed into
\begin{eqnarray*}
&&E\bigl[f\bigl(X_t^{x,\varepsilon}\bigr)\bigr] =E\bigl[f\bigl(
\bar{X}_t^{x,\varepsilon}\bigr)\bigr]+\sum
_{i=1}^m \varepsilon^i E \bigl[f\bigl(
\bar{X}_t^{x,\varepsilon}\bigr) \Phi_t^j
\bigr]+\varepsilon ^{m+1} R_m(t,x,\varepsilon),
\end{eqnarray*}
where:
\begin{longlist}[3.]
\item[1.]
%
\begin{eqnarray}\label{Rsmooth}
&&R_m(t,x,\varepsilon)\nonumber
\\
&&\qquad =\int_0^1(1-u)^m (m+1)
\\
&&\hspace*{45pt}{}\times \sum
_{\alpha^{(k)},\beta^{(k)}}^{m+1} E \Biggl[ \partial_{\alpha^{(k)}}
f\bigl(\tilde{X}_t^{x,\varepsilon u}\bigr) \prod
_{l=1}^k \frac{1}{\beta_l !} \frac{\partial^{\beta
_{l}}}{\partial\eta^{\beta_{l}}}
X_t^{x,\eta,\alpha_{l}}\bigg| _{\eta
=\varepsilon u} \Biggr]\,du
\nonumber
\end{eqnarray}
for $f \in C_b^{\infty}({\mathbf R}^N)$,
\item[2.]
\begin{eqnarray}\label{RLip}
\hspace*{-4pt}&&R_m(t,x,\varepsilon)
\nonumber
\\
\hspace*{-4pt}&&\qquad =\int_0^1(1-u)^m (m+1)
\nonumber\\[-8pt]\hspace*{-4pt} \\[-8pt]\nonumber
\hspace*{-4pt}&&\hspace*{45pt}{}\times  \sum_{\alpha^{(k)},\beta^{(k)}}^{m+1} E \Biggl[
\partial_{\alpha^{(1)}} f\bigl(\tilde{X}_t^{x,\varepsilon u}\bigr)
H_{\alpha
^{(k-1)}} \Biggl(Y_t^{\varepsilon u},\prod
_{l=1}^k \frac{1}{\beta_l
!} \frac{\partial^{\beta_{l}}}{\partial\eta^{\beta_{l}}}
X_t^{x,\eta,\alpha_{l}}\bigg| _{\eta=\varepsilon u} \Biggr) \Biggr]\,du
\nonumber
\end{eqnarray}
for $f \in C_b^{1}({\mathbf R}^N)$,
\item[3.]
%
\begin{eqnarray}\label{RBB}
\qquad &&R_m(t,x,\varepsilon)
\nonumber
\\
&&\qquad =\int_0^1(1-u)^m (m+1)
\\
&&\hspace*{45pt}{}\times \sum_{\alpha^{(k)},\beta^{(k)}}^{m+1} E \Biggl[ f\bigl(
\tilde{X}_t^{x,\varepsilon u}\bigr) H_{\alpha^{(k)}}
\Biggl(Y_t^{\varepsilon u},\prod_{l=1}^k
\frac{1}{\beta_l !} \frac
{\partial^{\beta_{l}}}{\partial\eta^{\beta_{l}}} X_t^{x,\eta,\alpha_{l}}\bigg|
_{\eta=\varepsilon u} \Biggr) \Biggr]\,du
\nonumber
\end{eqnarray}
for an arbitrary bounded continuous function $f$.
\end{longlist}

For $k\leq m+1$, $\sum_{l=1}^k \beta_l =m+1+k$, $\beta_l \geq2$,
$\alpha^{(k)} = (\alpha_1,\ldots, \alpha_k) \in\{1,\ldots,N \}
^k$, the product of the higher derivative terms with respect to
$\varepsilon$ of $X_t^{x,\varepsilon}$ is characterized as
%
\begin{equation}
\prod_{l=1}^k \frac{1}{\beta_l !}
\frac{\partial^{\beta
_{l}}}{\partial\varepsilon^{\beta_{l}}} X_t^{x,\varepsilon,\alpha
_{l}} \in{\mathcal
K}_{m+1+k}, \label{Korder}
\end{equation}
by using Lemma \ref{differentiationKr} with Lemma \ref{Properties}.

For $i=1,\ldots,N$ and $G(t,x) \in{\mathcal K}_r$, we are able to see
the following property for Malliavin weight as in Proposition~\ref{pr1}:
%
\begin{eqnarray}\label{MalliavinIBPY}
\qquad &&H_{(i)}\bigl(Y_t^{\varepsilon},G(t,x)\bigr)\nonumber
\\
&&\qquad =\delta
\Biggl(\sum_{j=1}^{N}G(t,x)
\gamma_{ij}^{Y_t^{\varepsilon
}}DY_t^{\varepsilon,j} \Biggr)
\\
&&\qquad = \Biggl[G(t,x) \sum_{j=1}^{N} \sum
_{k=1}^d \int_0^t
\gamma_{ij}^{Y_t^{\varepsilon}}\bigl(J_t^{x,\varepsilon
}\bigl(J_s^{x,\varepsilon}\bigr)^{-1}V_k
\bigl(X_s^{x,\varepsilon}\bigr)\bigr)^j
\,dB_s^k
\nonumber
\\
&&\hspace*{37pt}{} -\sum_{j=1}^{N}\sum
_{k=1}^d \int_0^t
\bigl[D_{s,k} G(t,x)\bigr] \gamma_{ij}^{Y_t^{\varepsilon
}}\bigl(J_t^{x,\varepsilon}\bigl(J_s^{x,\varepsilon
}\bigr)^{-1}V_k\bigl(X_s^{x,\varepsilon}\bigr)
\bigr)^j \,ds \Biggr] \nonumber
\\
&&\qquad \in{\mathcal K}_{r-1}.
\nonumber
\end{eqnarray}
Here, the first and the second terms in the second equality are
characterized by
%
\begin{eqnarray}
G(t,x) \sum_{j=1}^{N} \sum
_{k=1}^d \int_0^t
\gamma_{ij}^{Y_t^{\varepsilon}} \bigl(J_t^{x,\varepsilon
}
\bigl(J_s^{x,\varepsilon}\bigr)^{-1}V_k
\bigl(X_s^{x,\varepsilon}\bigr)\bigr)^j
\,dB_s^k &\in& {\mathcal K}_{r-1},
\\
\int_0^t \bigl[D_{s,k} G(t,x)
\bigr] \gamma_{ij}^{Y_t^{\varepsilon
}}\bigl(J_t^{x,\varepsilon}
\bigl(J_s^{x,\varepsilon
}\bigr)^{-1}V_k
\bigl(X_s^{x,\varepsilon}\bigr) \bigr)^j \,ds &\in& {\mathcal
K}_{r},
\end{eqnarray}
since
%
\begin{eqnarray}
&& \int_0^t \gamma_{ij}^{Y_t^{\varepsilon}}
\bigl(J_t^{x,\varepsilon
}\bigl(J_s^{x,\varepsilon}
\bigr)^{-1}V_k\bigl(X_s^{x,\varepsilon}\bigr)
\bigr)^j \,dB_s^k \in {\mathcal
K}_{-2+1}={\mathcal K}_{-1}.
\end{eqnarray}

Then, applying (\ref{MalliavinIBPY}) with
(\ref{Korder}) for (\ref{Rsmooth}), (\ref{RLip}) and (\ref{RBB}),
we obtain the following estimates
according to the smoothness of $f$:
\begin{longlist}[2.]
\item[1.]
%
\begin{equation}
\sup_{x\in\R^N} \bigl\llvert R_m(t,x,\varepsilon)
\bigr\rrvert \leq C \Biggl(\sum_{k=1}^{m+1}
t^{(m+1+k)/2}\bigl\llVert \nabla^k f \bigr\rrVert _{\infty}
\Biggr),
\end{equation}
for any $f \in C_b^{\infty}(\R^N)$,

\item[2.]
%
\begin{equation}
\sup_{x\in\R^N} \bigl\llvert R_m(t,x,\varepsilon)
\bigr\rrvert \leq C t^{(m+2)/2}\llVert \nabla f \rrVert _{\infty},
\end{equation}
for any $f \in C_b^1$,

\item[3.]
%
\begin{equation}
\sup_{x\in\R^N} \bigl\llvert R_m(t,x,\varepsilon)
\bigr\rrvert \leq C t^{(m+1)/2} \llVert f \rrVert _{\infty},
\end{equation}
for an arbitrary bounded continuous function $f$.
\end{longlist}
Then we have the assertion.

\section{Proof of Theorem~\texorpdfstring{\protect\ref{222222222222222222}}{2}}\label{appeB}
For $f \in C_b^{\infty}({\mathbf R}^N;{\mathbf R})$, we have
%
\begin{eqnarray}\label{formula1}
&&\int_{{\mathbf R}^N}f(y)E \biggl[H_{(i)} \biggl(
\frac{\partial
}{\partial\varepsilon} X_{t}^{x,\varepsilon}\bigg| _{\varepsilon=0}, G(t,x)
\biggr) \Big| \bar{X}_{t}^{x,\varepsilon}=y \biggr] \nu (dy)
\\[-1pt]
&&\qquad =E \biggl[ f \bigl( \bar{X}_{t}^{x,\varepsilon} \bigr)
H_{(i)} \biggl(\frac{\partial}{\partial\varepsilon} X_{t}^{x,\varepsilon}\bigg|
_{\varepsilon=0}, G(t,x) \biggr) \biggr]
\nonumber
\\[-1pt]
&&\qquad =E \Biggl[ f \bigl( \bar{X}_{t}^{x,\varepsilon} \bigr) \delta
\Biggl(\sum_{j=1}^{N}G(t,x)
\gamma_{ij}^{Y_t^{0}}DY_t^{0,j} \Biggr)
\Biggr]
\nonumber
\\[-1pt]
&&\qquad =E \Biggl[ f \bigl( \bar{X}_{t}^{x,\varepsilon} \bigr) \delta
\Biggl(\sum_{j=1}^{N}\varepsilon G(t,x)
\frac{1}{\varepsilon^2} \gamma_{ij}^{Y_t^{0}} \varepsilon
DY_t^{0,j} \Biggr) \Biggr]
\nonumber
\\[-1pt]
&&\qquad = E \Biggl[ f \bigl( \bar{X}_{t}^{x,\varepsilon} \bigr) \delta
\Biggl(\sum_{j=1}^{N}\varepsilon G(t,x)
\gamma_{ij}^{\bar
{X}_t^{x,\varepsilon}} D \bar{X}_t^{x,\varepsilon,j}
\Biggr) \Biggr]
\nonumber
\\[-1pt]
&&\qquad = E \bigl[ f \bigl( \bar{X}_{t}^{x,\varepsilon}
\bigr)H_{(i)} \bigl( \bar{X}_{t}^{x,\varepsilon}, \varepsilon
G(t,x) \bigr) \bigr]
\nonumber
\\[-1pt]
&&\qquad = E \bigl[ \partial_i f \bigl( \bar{X}_{t}^{x,\varepsilon}
\bigr) \varepsilon G(t,x) \bigr]
\nonumber
\\[-1pt]
&&\qquad = \int_{{\mathbf R}^N}\partial_i f(y)E \bigl[
\varepsilon G(t,x) | \bar{X}_{t}^{x,\varepsilon}=y \bigr] \nu(dy)
\nonumber
\\[-1pt]
&&\qquad =\int_{{\mathbf R}^N} f(y) \partial_i^{\ast}E
\bigl[\varepsilon G(t,x) | \bar{X}_{t}^{x,\varepsilon}=y \bigr]
\nu(dy), \label{formula2}
\end{eqnarray}
where $\gamma^{\bar{X}_t^{x,\varepsilon}}=(\gamma_{ij}^{\bar
{X}_t^{x,\varepsilon}} )_{1\leq i,j \leq N}$ and $\gamma^{Y_t^{0}}=(
\gamma_{ij}^{Y_t^{0}} )_{1\leq i,j \leq N}$ are the inverse matrices
of the Malliavin covariance matrices of $\bar{X}_t^{x,\varepsilon}$
and ${Y}_t^{0}$, respectively.
Here, we note that $Y_t^0=\frac{\partial}{\partial\varepsilon}
X_{t}^{x,\varepsilon}\mid _{\varepsilon=0}$ and
$\bar{X}_{t}^{x,\varepsilon}=X_t^{x,0}+\varepsilon\frac{\partial
}{\partial\varepsilon} X_{t}^{x,\varepsilon}\mid _{\varepsilon
=0}=X_t^{x,0}+\varepsilon Y_t^0$.
Also, we use the following relations in the above equations; for
$k=1,\ldots,d$ and $j=1,\ldots,N$,
%
\begin{equation}
D_{s,k} \bar{X}_t^{x,\varepsilon,j}=\varepsilon
D_{s,k} \frac
{\partial}{\partial\varepsilon} X_t^{x,\varepsilon,j}\bigg|
_{\varepsilon=0}= \varepsilon D_{s,k} Y_t^{0,j},\qquad
s \leq t,
\end{equation}
and, for $i,j=1,\ldots,N$,
%
\begin{equation}
\gamma_{ij}^{\bar{X}_t^{x,\varepsilon}}=\frac{1}{\varepsilon^2} \gamma_{ij}^{Y_t^{0}}.
\end{equation}
Formulas (\ref{formula1}) and (\ref{formula2}) hold for any
Lipschitz and bounded Borel function $f$ by using mollifier arguments.
We remark that in general for any $G \in{\mathbf D}^{\infty}$ and
nondegenerate $F \in{\mathbf D}^{\infty}({\mathbf R}^N)$, the conditional
expectation can be regarded as a map ${\mathbf D}^{\infty} \ni G \mapsto
E[G\mid F=\cdot ] \in{\mathcal S}({\mathbf R}^N)$ by \citet{Malliavin} and
\citet{Malliavin-Thalmaier}. Therefore, for
$k=1,\ldots,j \leq m$, $\sum_{l=1}^k \beta_l =j+k$, $\beta_l \geq
2$, $\alpha^{(k)} = (\alpha_1,\ldots, \alpha_k) \in\{1,\ldots,N \}^k$,
we have
\begin{eqnarray}
&&E \Biggl[H_{{\alpha}^{(k)}} \Biggl(\frac{\partial}{\partial
\varepsilon} X_{t}^{x,\varepsilon}
\bigg| _{\varepsilon=0}, \prod_{l=1}^k
\frac{1}{\beta_l !}\frac{\partial^{\beta_{l}}}{\partial
\varepsilon^{\beta_{l}}} X_t^{x,\varepsilon,\alpha
_{l}}\bigg|
_{\varepsilon=0} \Biggr) \Big| \bar{X}_{t}^{x,\varepsilon
}=y \Biggr]
\nonumber\\[-8pt]\\[-8pt]\nonumber
&&\qquad =  \varepsilon^k \partial^{\ast}_{\alpha_k} \circ
\partial^{\ast
}_{\alpha_{k-1}} \circ\cdots\circ\partial^{\ast}_{\alpha_1}
E \Biggl[ \prod_{l=1}^k
\frac{1}{\beta_l !}\frac{\partial^{\beta
_{l}}}{\partial\varepsilon^{\beta_{l}}} X_t^{x,\varepsilon,\alpha
_{l}}\bigg|
_{\varepsilon=0} \Big| \bar{X}_{t}^{x,\varepsilon} =y \Biggr],
\nonumber
\end{eqnarray}
and obtain the assertion.

\section{Proof of Theorem~\texorpdfstring{\protect\ref{33333333333333}}{3}}\label{appeC}
We follow the similar argument as in \citeauthor{Kusuokaa} (\citeyear{Kusuokaa,Kusuokac,Kusuokad}) and
Chapter~3 of \citet{Crisan-Manolarakis-Nee}.

Note first that we have the following equality:
\begin{eqnarray*}
&&P_Tf(x) - Q^m_{(s_n)} Q^m_{(s_{n-1})}
\cdots Q^m_{(s_1)} f(x)
\\[-2pt]
&&\qquad =P_{T-t_{n-1}}P_{t_{n-1}} f(x) - Q^m_{(s_n)}
P_{t_{n-1}} f(x)
\\[-2pt]
&&\quad\qquad{}+Q^m_{(s_n)}P_{t_{n-1}} f(x)-Q^m_{(s_n)}Q^m_{(s_{n-1})}
P_{t_{n-2}}f (x)
\\[-2pt]
&&\quad\qquad{}+\cdots
\\[-3pt]
&&\quad\qquad{}+Q^m_{(s_{n})}\cdots Q^m_{(s_2)}P_{t_1}
f-Q^m_{(s_n)} Q^m_{(s_{n-1})}\cdots
Q^m_{(s_1)}f
\\[-2pt]
&&\qquad =P_{T-t_{n-1}}P_{t_{n-1}} f(x) - Q^m_{(s_n)}
P_{t_{n-1}} f(x)
\\[-2pt]
&&\quad\qquad {}+Q^m_{(s_n)} \bigl(P_{s_{n-1}}P_{t_{n-2}}
f(x)-Q^m_{(s_{n-1})} P_{t_{n-2}}f (x)\bigr)
\\[-4pt]
&&\quad\qquad{}+ \cdots
\\[-3pt]
&&\quad\qquad{}+Q^m_{(s_{n})}\cdots Q^m_{(s_2)}
\bigl(P_{t_1} f(x)-Q^m_{(s_1)}f(x)\bigr).
\end{eqnarray*}
Then, since $Q^m$ is a Markov operator, we have
\begin{eqnarray*}
&&\bigl\llVert P_Tf - Q^m_{(s_n)}
Q^m_{(s_{n-1})}\cdots Q^m_{(s_1)} f\bigr
\rrVert _{\infty}
\\[-2pt]
&&\qquad \leq\bigl(\bigl\llVert P_{s_{n}}P_{t_{n-1}} f - Q^m_{(s_n)}
P_{t_{n-1}} f \bigr\rrVert _{\infty
}
\\[-2pt]
&&\hspace*{3pt}\quad\qquad{} + \bigl\llVert P_{s_{n-1}}P_{t_{n-2}} f -Q^m_{(s_{n-1})}
P_{t_{n-2}}f \bigr\rrVert _{\infty}
\\[-4pt]
&&\hspace*{5pt}\quad\qquad\cdots
\\[-4pt]
&&\hspace*{66pt}\quad\qquad{}+ \bigl\llVert P_{t_1} f-Q^m_{(s_1)}f\bigr
\rrVert _{\infty}\bigr)\bigl(1+\mathrm{O}(\varepsilon)\bigr)
\\[-2pt]
&&\qquad =  \Biggl(\sum_{k=2}^n \bigl\llVert
P_{s_{k}}P_{t_{k-1}} f -Q^m_{(s_{k})}
P_{t_{k-1}}f \bigr\rrVert _{\infty}
\\[-2pt]
&&\hspace*{54pt}\quad\qquad{}+ \bigl\llVert P_{t_1} f-Q^m_{(s_1)}f\bigr
\rrVert _{\infty}\Biggr)\bigl(1+\mathrm{O}(\varepsilon)\bigr).
\end{eqnarray*}

First, note that we can directly apply (\ref{shorttimeerror1}), (\ref
{shorttimeerror2}) or (\ref{shorttimeerror3}) in
Corollary~\ref{shorttimeerror} to obtain an estimate of
$\llVert   P_{t_1} f-Q^m_{(s_1)}f\rrVert  _{\infty}$ for $f \in C_b^{\infty}(\R
^N;{\mathbf R})$,
a Lipschitz\vspace*{2pt} continuous function or a bounded Borel function, respectively.
To obtain an estimate of
$\sum_{k=2}^n \llVert   P_{s_{k}}P_{t_{k-1}} f -Q^m_{(s_{k})} P_{t_{k-1}}f \rrVert
_{\infty}$,
we apply the results in Corollary~\ref{shorttimeerror} to $P_t f$ (in
stead of $f$) as follows:
\begin{itemize}
\item
By (\ref{shorttimeerror1}) in Corollary~\ref{shorttimeerror},
for $s,t \in(0,1]$ and $f \in C_b^{\infty}(\R^N;{\mathbf R})$, there
exists $C$ such that
%
\begin{eqnarray}
\bigl\llVert P_{s}P_{t} f -Q^m_{(s)}
P_{t}f \bigr\rrVert _{\infty} &\leq& \sum
_{l=1}^{m+1} s^{(m+1+l)/2} C \bigl\llVert
\nabla^l P_t f \bigr\rrVert _{\infty
}
\\[-2pt]
&\leq& \sum_{l=1}^{m+1} s^{(m+1+l)/2} C
\bigl\llVert \nabla^l f \bigr\rrVert _{\infty}.
\end{eqnarray}
Hence,
%
\begin{eqnarray}
&&\bigl\llVert P_Tf - Q^m_{(s_n)}
Q^m_{(s_{n-1})}\cdots Q^m_{(s_1)} f\bigr
\rrVert _{\infty
}
\\
&&\qquad \leq C \sum_{k=2}^n \sum
_{l=1}^{m+1} s_k^{(m+1+l)/2} \bigl
\llVert \nabla^l f \bigr\rrVert _{\infty}\label{21}
\\
&&\quad\qquad{}+ C \sum_{l=1}^{m+1} s_1^{(m+1+l)/2}
\bigl\llVert \nabla^l f \bigr\rrVert _{\infty}.
\end{eqnarray}

\item
By (\ref{shorttimeerror2}) in Corollary \ref{shorttimeerror},
for $s,t \in(0,1]$ and $f \in C_b^{1}(\R^N;{\mathbf R})$, there exists
$C$ such that
%
\begin{eqnarray}
\bigl\llVert P_{s}P_{t} f -Q^m_{(s)}
P_{t}f \bigr\rrVert _{\infty} &\leq& s^{(m+2)/2} C \llVert
\nabla P_t f \rrVert _{\infty}
\\
&\leq& s^{(m+2)/2} C \llVert \nabla f \rrVert _{\infty}.
\end{eqnarray}
Hence,
%
\begin{eqnarray}
&&\bigl\llVert P_Tf - Q^m_{(s_n)}
Q^m_{(s_{n-1})}\cdots Q^m_{(s_1)} f\bigr
\rrVert _{\infty
}
\\
&&\qquad \leq  C \sum_{k=2}^n
s_k^{(m+2)/2} \llVert \nabla f \rrVert _{\infty}\label
{22}
\\
&&\quad\qquad{}+C s_1^{(m+2)/2} \llVert \nabla f \rrVert _{\infty}.
\end{eqnarray}

\item
By (\ref{shorttimeerror3}) in Corollary \ref{shorttimeerror},
for $s,t \in(0,1]$ and bounded Borel function $f$ on $\R^N$, there
exists $C$ such that
%
\begin{eqnarray}
\bigl\llVert P_{s}P_{t} f -Q^m_{(s)}
P_{t}f \bigr\rrVert _{\infty} &\leq& s^{(m+1)/2} C \llVert
P_t f \rrVert _{\infty}
\\
&\leq&s^{(m+1)/2} C \llVert f \rrVert _{\infty}.
\end{eqnarray}
Hence,
%
\begin{eqnarray}
&&\bigl\llVert P_Tf - Q^m_{(s_n)}
Q^m_{(s_{n-1})}\cdots Q^m_{(s_1)} f\bigr
\rrVert _{\infty
}
\\
&&\qquad \leq C \sum_{k=2}^n
s_k^{(m+1)/2} \llVert f \rrVert _{\infty}\label{23}
\\
&&\quad\qquad{}+C s_1^{(m+1)/2} \llVert f \rrVert _{\infty}.
\end{eqnarray}
\end{itemize}

Next, we obtain more explicit and compact expressions with regard to $n$
particularly for (\ref{21}), (\ref{22}) and (\ref{23}).

First, from the definition of $s_k$ for $k \in\{2,\ldots,n \}$, we have
%
\begin{equation}
s_k = \frac{\gamma T (k-1)^{\gamma-1}}{n^{\gamma}} \int_{k-1}^k
\bigl(u/(k-1)\bigr)^{\gamma-1}\,du.
\end{equation}
For $k \in\{2,\ldots,n \}$, $(u/(k-1))^{\gamma-1}\leq\max\{
(k/(k-1))^{\gamma-1},1 \} \leq\max\{2^{\gamma-1},1 \}$. Then
%
\begin{eqnarray}
s_k^{l/2} &\leq& \biggl(\frac{\gamma T (k-1)^{\gamma-1}}{n^{\gamma
}}\max\bigl
\{2^{\gamma-1},1 \bigr\} \biggr)^{l/2}
\\
&\leq& C (1/n)^{\gamma l/2} (k-1)^{(\gamma-1) l /2},
\end{eqnarray}
where $C=C(T,\gamma)$.

We consider the estimates for three different ranges of $\gamma$
that are larger than, equal to and less than $(l-2)/l$, respectively.
[$\gamma= (l-2)/l$ satisfies $(\gamma-1)l/2=-1$.]

For $0< \gamma< (l-2)/l$,
%
\begin{equation}
C (1/n)^{\gamma l/2} \sum_{k=2}^n
(k-1)^{(\gamma-1)l/2} \leq C (1/n)^{\gamma l/2}.
\end{equation}
For $\gamma=(l-2)/l$,
%
\begin{eqnarray}
&&C (1/n)^{\gamma l/2} \sum_{k=2}^n
(k-1)^{(\gamma-1)l/2}
\\
&&\qquad =C(1/n)^{(l-2)/2} \sum_{k=1}^n
(k-1)^{-1}
\\
&&\qquad \leq C (1/n)^{(l-2)/2} \log n.
\end{eqnarray}
For $\gamma>(l-2)/l$,
%
\begin{eqnarray}
&&C (1/n)^{\gamma l/2} \sum_{k=2}^n
(k-1)^{(\gamma-1)l/2}
\\
&&\qquad =C (1/n)^{(\gamma-1)l/2}(1/n)^{l/2} \sum
_{k=2}^n (k-1)^{(\gamma
-1)l/2}
\\
&&\qquad =C (1/n)^{(l-2)/2} \sum_{k=2}^n
\biggl(\frac{k-1}{n} \biggr)^{(\gamma-1)l/2} \frac{1}{n}
\\
&&\qquad \leq C (1/n)^{(l-2)/2}.
\end{eqnarray}
Then,\vspace*{1pt}
by combining an estimate of
$\llVert   P_{t_1} f-Q^m_{(s_1)}f\rrVert  _{\infty}$ for $f \in C_b^{\infty}(\R
^N;{\mathbf R})$,
a Lipschitz continuous function or a bounded Borel function,
we have the assertion.
\end{appendix}

\section*{Acknowledgements}
We are very grateful to the Editor, the Associate Editor
and two anonymous referees for their precious comments and suggestions.


%

\printaddresses
\end{document}